%% file: KKCC3Drevision.tex
\documentclass{article}[12pt]  
\usepackage[utf8x]{inputenc}
\usepackage{amsfonts}
\usepackage{graphics}
\usepackage{graphicx} 
\usepackage{color}
\usepackage[margin=3.cm]{geometry}
\usepackage{amsmath,amssymb}
\usepackage{color}
\usepackage[dvipsnames]{xcolor}
\usepackage{appendix}
\usepackage{amsthm}
\usepackage{array}
\usepackage{mathtools}
\renewcommand{\epsilon}{\varepsilon}

\graphicspath{{figures/}}

\DeclareGraphicsExtensions{.eps,.pdf,.png,.jpg,.JPG,.jpeg,.ps,.pdf_tex}
\numberwithin{equation}{section}

\usepackage{setspace,color}
\usepackage{amsmath,amssymb,mathrsfs}
\usepackage{graphicx,epstopdf}
\usepackage{xcolor}

\numberwithin{equation}{section}

\newcommand{\ihat}{\hat{\i}}
\newcommand{\jhat}{\hat{\j}}
\newcommand{\khat}{\hat{{\rm k}}}

\begin{document}

 \title{Close evaluation of layer potentials in three dimensions}
\author{Shilpa Khatri\footnote{Applied Mathematics Unit, School of Natural Sciences,
    University of California, Merced, 5200 North Lake Road, Merced, CA
    95343} \and  Arnold D. Kim ${}^{\ast}$ \and R. Cortez\footnote{Mathematics Department, Tulane University, 424 Gibson Hall, New Orleans, LA 70118} \and Camille Carvalho ${}^{\ast}$}
   \maketitle

  \begin{abstract}
    We present a simple and effective method for evaluating double-
    and single-layer potentials for Laplace's equation in three
    dimensions close to the boundary. The close evaluation of these
    layer potentials is challenging because they are nearly singular
    integrals. The method we propose is based on writing these layer
    potentials in spherical coordinates where the
    point at which their kernels are peaked maps to the north pole. An
    $N$-point Gauss-Legendre quadrature rule is used for integration
    with respect to the the polar angle rather than the cosine of the
    polar angle. A $2N$-point periodic trapezoid rule is used to
    compute the integral with respect to the azimuthal angle which
    {acts as} a natural and effective averaging operation in this
    coordinate system. The numerical method resulting from combining
    these two quadrature rules in this rotated coordinate system
    yields results that are consistent with asymptotic behaviors of
    the double- and single-layer potentials at close evaluation
    distances. In particular, we show that the error in computing the
    double-layer potential{, after applying a subtraction method,} is
    quadratic with respect to the evaluation distance from the
    boundary, and the error is linear for the single-layer
    potential. We improve upon the single-layer potential by
    introducing an alternate approximation based on a perturbation
    expansion and obtain an error that is quadratic with respect to
    the evaluation distance from the boundary.
  \end{abstract}

  \textbf{Keywords: }
    Nearly singular integrals, close evaluation problem, potential
    theory, boundary integral equations, numerical quadrature.

\section{Introduction}

The close evaluation problem arises when using boundary integral methods
{to solve} boundary value problems for linear, elliptic partial
differential equations. In boundary integral methods, the solution of
the boundary value problem is given in terms of double- and
single-layer potentials, integrals of a kernel multiplied by a density
over the boundary of the domain. The kernel for the single-layer
potential is the fundamental solution of the elliptic partial
differential equation and the kernel for the double-layer potential is
the normal derivative of that fundamental solution. Each of these
kernels has an isolated singularity at a known point on the boundary.
When evaluating layer potentials at points close to the boundary, the
associated kernel is regular, but is sharply peaked. For
this reason, we say that layer potentials evaluated close to the
boundary are nearly singular integrals.

Nearly singular integrals are more challenging to compute numerically
than weakly singular ones, resulting from evaluating layer potentials at points on the boundary.  When computing weakly singular integrals, one explicitly
addresses the singularity in the kernel. There are several high-order
methods available to compute weakly singular integrals
(e.g. \cite{atkinson1997numerical,bruno2001fast,graham2002fully,ganesh2004high,bremer2010nonlinear}).
For example,
  {they} can be computed accurately
using high-order product Nystr\"om
methods~\cite{atkinson1997numerical, delves1988computational} that
analytically treat integration over the singular point. These product
Nystr\"om methods are often used to solve the boundary integral
equations for the density. However, for a nearly singular integral,
there is no singularity to address. Nonetheless, without effectively
addressing the peaked behavior in a nearly singular integral for an
evaluation point fixed close to the boundary, one must increase the
number of quadrature points to obtain accuracy commensurate with
evaluation points that are far from the boundary.

There are two factors that affect the accuracy of a numerical method
for the close evaluation problem: the distance from the boundary
 and the number of quadrature points. Based on these two
factors, one can design a new competitive numerical method (improving
the convergence rate with respect to the number of quadrature points), or one can provide
corrections to existing ones (improving the convergence rate with
respect to the distance for a fixed number of quadrature points).  Here, we focus on this second
consideration.

The close evaluation problem has been studied extensively for
two-dimensional problems.  Schwab and
Wendland~\cite{schwab1999extraction} have developed a boundary
extraction method based on a Taylor series expansion of the layer
potentials.  Beale and Lai~\cite{beale2001method} have developed a
method that first regularizes the nearly singular kernel and then adds
corrections for both the discretization and the regularization. The
result of this approach is a uniform error in space.  Helsing and
Ojala~\cite{helsing2008evaluation} developed a method that combines a
globally compensated quadrature rule and interpolation to achieve very
accurate results over all regions of the
domain. Barnett~\cite{barnett2014evaluation} has used surrogate local
expansions with centers placed near, but not on, the boundary. These results 
led to {the} work by Kl\"{o}ckner {\it et
  al.}~\cite{klockner2013quadrature} that introduces Quadrature By
Expansion (QBX). QBX uses expansions about accurate evaluation points
far away from the boundary to compute accurate evaluations for points close to
it. The convergence of QBX has been studied
in~\cite{epstein2013convergence}. Moreover, fast implementations of
QBX have since been developed~\cite{af2016fast, rachh2017fast,
  wala20183dqbx}, and rigorous error estimates have been derived for
the method~\cite{af2017error}. Recently, the authors have developed a
method that involves matched asymptotic expansions for the
kernel {of layer potentials}~\cite{carvalho2018asymptotic}. In that method, the asymptotic
expansion that captures the peaked behavior of the kernel can be
integrated exactly using Fourier series, and the relatively smooth
remainder is integrated numerically, resulting in a highly accurate
method.

There are fewer results for three-dimensional problems.  Beale {\it et
  al.}~\cite{beale2016simple} have extended the regularization method
to three-dimensional problems.  Additionally, QBX has been used for
three-dimensional problems~\cite{af2016fast, wala20183dqbx}. In
principle, the matched asymptotic expansion method developed by the authors  for two dimensions can
be extended to three-dimensional problems (using spherical harmonics
instead). However, we do not pursue that approach because the method
we present here is 
direct and
 {simple} to implement. 
 {The development of accurate numerical methods in three dimensions can be challenging and there is a need for such methods that are straightforward to implement.}

{Further, there exist quadrature methods to specifically deal with nearly singular integrals. Johnston and Elliot~\cite{johnston2005sinh} introduce a hyperbolic sine transformation to cluster quadrature points 
around the near singularity of the kernel. This clustering of points is dependent on the peakedness of the kernel and location of the singularity. Iri, Moriguti, and Takasawa~~\cite{iri1987certain} have developed a method for nearly singular integrals that takes advantage of the spectral accuracy of the periodic trapezoid rule. This method also results in a clustering of quadrature points around the near singularity of the kernel.
 }

In this paper, we study the close evaluation of double- and
single-layer potentials for Laplace's equation in three
dimensions. 
%
%
  We compute the double- and single-layer potentials in spherical
  coordinates where the isolated singular point on the boundary is
  aligned with the north pole{, possibly after a rotation}. At close evaluation distances, we show
  that the leading asymptotic behaviors for the kernels of both the
  double- and single-layer potentials are azimuthally invariant. It
  follows that integrating over the azimuth in this spherical coordinate
  system introduces a natural averaging operation that effectively
  smooths the nearly singular behavior of the layer potentials. We use
  these asymptotic results to introduce a simple and effective numerical method
  for computing the close evaluation of double- and single-layer
  potentials. {Comparisons with other quadrature methods designed for nearly singular integrals show that the proposed method is more accurate.}
%

%
%

We precisely define the close evaluation problem for the double and
single-layer potentials in Section \ref{sec:close-eval}. {We present in Section \ref{sec:prior} some prior quadrature methods for the close evaluation problem}. In Section
\ref{sec:local-analysis}, we derive the asymptotic behavior of the
contributions to {the double and single-}layer potentials made by a small region of the
boundary containing the point where the kernels are peaked.  By doing
so, we find a natural numerical method to evaluate these nearly
singular integrals. This numerical method is given in Section
\ref{sec:numerics}{, and improvements for the single-layer potential are considered in Section \ref{sec:extension}}. Several examples demonstrating the accuracy of
this numerical method are presented in Section
\ref{sec:examples}. Section \ref{sec:conclusions} gives our
conclusions. The details of the method we use to rotate the
integrals in spherical coordinates is given in \ref{sec:rotation}.

\section{Close evaluation of layer potentials}
\label{sec:close-eval}

Consider a simply connected, open set denoted by
$D \subset \mathbb{R}^{3}$ with {smooth oriented} boundary, $B$, and let
$\bar{D} = D \cup B$.  The function
$u \in C^{2}(D) \cap C^{1}(\bar{D})$ satisfying Laplace's equation,
$\Delta u = 0$, can be expressed using a \textit{representation
  formula}, a combination of double- and single-layer
potentials~\cite{guenther1996partial}:

\begin{equation}
  u(x) = - \frac{1}{4 \pi} \int_{B} \frac{n(y) \cdot ( x
    - y )}{|x - y|^{3}} u(y) \, \mathrm{d}\sigma_{y} + \frac{1}{4\pi}
  \int_{B} \frac{1}{|x - y|} \partial_{n} u(y)
\,  \mathrm{d}\sigma_{y},
  \quad x \in D,
  \label{eq:repres_int}
\end{equation}
where $n(y)$ is the outward unit normal at $y \in B$ and
$\mathrm{d}\sigma_{y}$ is the 
  surface element.

Throughout this paper, we consider \eqref{eq:repres_int}, which 
corresponds to the representation formula for interior problems.  To
consider the representation formula for exterior problems,  $n(y)$ can 
just be replaced by $-n(y)$ in \eqref{eq:repres_int}.  The solution of
the interior Dirichlet problem for Laplace's equation is often
represented using only the double-layer potential (first term in
\eqref{eq:repres_int}) with $u$ replaced by a general density
$\mu$. The solution of the exterior Neumann problem for Laplace's
equation is often represented using only the single-layer potential
(second term in \eqref{eq:repres_int}) with $\partial_{n} u$ replaced
by a general density $\rho$. For those problems, the densities $\mu$
and $\rho$ satisfy boundary integral equations which can be solved
accurately using Nystr\"om methods~\cite{atkinson1997numerical,
  delves1988computational}. Nonetheless, the close evaluation problem
for the double-layer potential for the interior Dirichlet problem and for
the single-layer potential for the exterior Neumann problem still
remains {a subject of active investigation}.  Here, we choose to study the close evaluation of
\eqref{eq:repres_int} since it includes both the close evaluation of
the double- and single-layer potentials.

\begin{figure}[h!]
  \centering
  \def\svgwidth{0.3\columnwidth} 
  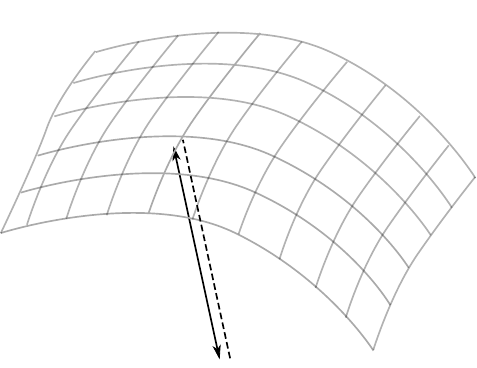
  \caption{Schematic presenting the notation given in \eqref{eq:closept}.}
  \label{fig:close-eval}
\end{figure}

Assuming that both $u$ and $\partial_{n} u$ on $B$ are known, one can
use high-order quadrature rules~\cite{delves1988computational} to
evaluate \eqref{eq:repres_int} with high accuracy. However, this high
accuracy is lost for evaluation points that are close to, but off the boundary. To
understand why this happens, we write a close evaluation point
$x \in D$ as,
\begin{equation}
  x = y^{\star} - \epsilon \ell n^{\star},
  \label{eq:closept}
\end{equation}
where $0 < \epsilon \ll 1$ is a small, dimensionless
parameter, $y^{\star} \in B$ is the point closest to $x$ on the
boundary, $n^{\star} = n(y^{\star})$ is the unit, outward normal from
$D$ at $y^{\star}$, and $\ell$ is a characteristic length of $B$ {associated with $y^*$, such as the radius of mean
curvature},  as shown in Fig. \ref{fig:close-eval}. 
  
  By substituting
\eqref{eq:closept} into \eqref{eq:repres_int}, and making use of
Gauss' law~\cite{guenther1996partial},
\begin{equation}
  \frac{1}{4 \pi} \int_{B} \frac{n(y) \cdot ( x - y
    )}{|x - y|^{3}} \mathrm{d}\sigma_{y} = \begin{cases}
    -1 & x \in D,\\
    -\frac{1}{2} & x \in B,\\
    \,\,\,\, 0 & x \in \mathbb{R}^{3} \backslash \bar{D}, \end{cases}
  \label{eq:gauss}
\end{equation}
we can rewrite the representation formula, \eqref{eq:repres_int}, evaluated at the close
evaluation point, \eqref{eq:closept}, as
\begin{equation}
  u(y^{\star} - \epsilon \ell n^{\star}) = 
  u(y^{\star}) - \frac{1}{4 \pi} \int_{B} \frac{n(y) \cdot
    ( y^{\star} - y - \epsilon \ell n^{\star})}{| y^{\star} - y -
    \epsilon \ell n^{\star} |^{3}} \left[ u(y) - u(y^{\star})
  \right] \mathrm{d}\sigma_{y} + \frac{1}{4\pi} \int_{B}
  \frac{1}{|y^{\star} - y - \epsilon \ell n^{\star}|} \partial_n u (y) \,
  \mathrm{d}\sigma_{y}.
  \label{eq:repres-int-sub}
\end{equation}
We call \eqref{eq:repres-int-sub} the modified representation formula
resulting from applying a subtraction
method~\cite{delves1988computational} to the double-layer potential in \eqref{eq:repres_int}.

The close evaluation of \eqref{eq:repres-int-sub} corresponds to 
the asymptotic limit, $\epsilon \to 0^{+}$ and in that case the integrals in
\eqref{eq:repres-int-sub} are nearly singular. Setting $\epsilon = 0$ in
\eqref{eq:repres-int-sub} we find that the kernels become singular at
$y = y^{\star}$, but those singularities are integrable. 
Nearly singular
integrals place high demands on any numerical integration method 
since there are no
explicit singularities that can be treated
analytically. Consequently, 
these nearly singular integrals require the development of 
{specialized}
methods to handle the sharply peaked integrands, which is challenging if the error is to be uniformly bounded over the entire domain. 

{ In our analysis and presentation of a  new method,
we assume $B$ is an analytic and closed oriented surface that can be parameterized using spherical coordinates, $y = y(s,t)$ with $s \in [0, \pi]$ and $t \in [-\pi, \pi]$, where $s$ denotes the polar angle, and $t$ denotes the azimuthal angle. In these coordinates, we assume we have rotated the domain to set $y^\star = y(0, \cdot)$, where $y(0,\cdot)$
denotes the spherical mean, $y(0,\cdot):= \lim \limits_{s \to 0^+}
\displaystyle\frac{1}{2\pi} \int_{-\pi}^{\pi} y(s,t) \,
dt$. In other words, $y^\star$ is aligned with the north pole, where we make use of the spherical mean to define {all} quantities. {This notation will be used throughout the paper to denote the spherical mean.}  Note that one can always apply a rotation so that any chosen $y^\star$ maps to $y(0, \cdot)$. The expression of the transformation matrix we use for this rotation can be found in \ref{sec:rotation}.
Using this parametrization and rotation, we write the integrals in (\ref{eq:repres-int-sub}) as, 
\begin{equation}
  I(y^{\star}) = \frac{1}{4\pi} \int_{-\pi}^{\pi}
  \int_{0}^{\pi} \tilde{F}(s,t)
 \, \mathrm{d}s \mathrm{d}t \, ,
 \label{eq:integrals}
\end{equation}
with
\begin{multline}
  \tilde{F}(s,t ) = -
   \frac{\tilde{n}(s,t) \cdot ( y(0,\cdot) - y(s,t) - \epsilon \ell
   n^\star)}{| y(0,\cdot) - y(s,t) - \epsilon \ell
      n^\star|^{3}}
  \tilde{J}(s,t) \left[ \tilde{\mu}(s,t) - \tilde{\mu}(0,\cdot) \right]
  + \frac{1}{|y(0,\cdot) - y(s,t) - \epsilon \ell
     n^\star|} \tilde{J}(s,t) \tilde{\rho}(s,t) \, ,
\end{multline}
where $\tilde{J}(s,t) = | \partial_{s} y(s,t) \times \partial_{t} y(s,t)
|$, $\tilde{n}(s,t) = n(y(s,t))$, $\tilde{\mu}(s,t) = \mu(y(s,t))$, and $\tilde{\rho}(s,t) = \rho(y(s,t))$. Note
that $n^\star = \tilde{n}(0,\cdot) $. 
For simplicity, we have replaced $u$ by $\mu$ and $\partial_n u$ by $\rho$, and we assume $\mu$ and $\rho$ are known and smooth.
}

{
\section{Prior quadrature methods developed for the close evaluation problem} 
\label{sec:prior} 
In this manuscript, we present  a new numerical method for computing \eqref{eq:repres-int-sub}
that  is based on an asymptotic analysis as $\epsilon \to 0^{+}$. 
As stated in the introduction, quadrature methods have been previously developed to accurately approximate nearly singular
integrals. In this section, we give three examples of commonly used existing methods that we will compare to our new method in Section \ref{sec:examples}. }

{
The three quadrature rules we present here are (1) the product Gauss quadrature rule by Atkinson~\cite{atkinson1982numerical} (PGQ) which uses Gauss-Legendre quadrature abscissas in the polar angle, (2) a hyperbolic sine tranformation presented by Johnston and Elliot in~\cite{johnston2005sinh} (SINH), and (3) a quadrature formula presented by Iri, Moriguti, and Takasawa in~\cite{iri1987certain} with the implementation as given in Robinson and Doncker~\cite{robinson1981algorithm} (IMT). All three methods are defined on a domain, $(z,t) \in [-1,1] \times [-\pi, \pi]$, and use the common substitution for the polar angle for integrals in spherical coordinates, $z = \cos(s)$ where $s \in [0, \pi]$~\cite{atkinson1982numerical}. 
Let
$z_{i} \in (-1,1)$ for $i = 1, \cdots, N$ denote the $N$-point quadrature rule abscissas with corresponding weights,
$w_{i}$ and then $s_i = \cos^{-1}(z_{i})$ for the abscissas in the polar angle.  
For all three methods, we use the periodic trapezoid rule in the azimuthal direction. Let $t_{j} = -\pi + \pi (j - 1)/N \in  [-\pi,\pi]$ for
$j = 1, \cdots, 2N$ denote the equi-spaced grid points for the periodic trapezoid rule. This results in the numerical method 
\begin{equation}
  I(y^{\star}) \approx I^{N}(y^{\star}) = \frac{1}{4 N} \sum_{i =
    1}^{N} \sum_{j = 1}^{2N} w_{i} 
\tilde{F}(s_{i},t_{j}) \, ,
  \label{eq:genquadrature}
\end{equation}
to evaluate (\ref{eq:integrals}) where $s_i$ and $w_i$ differ between the three methods. }

{ We present numerical results for the solution of Laplace's equation using the modified representation formula, \eqref{eq:repres-int-sub}, for close evaluation points interior to a peanut-shaped domain presented in Atkinson~\cite{atkinson1982laplace,
  atkinson1985algorithm}. The boundary of this domain is given by, 
\begin{equation}
  y(\theta,\varphi) = r(\theta) ( \sin\theta \cos\varphi, 2
  \sin\theta \sin\varphi, \cos\theta ), \quad \theta \in [0,\pi],
  \quad \varphi \in [-\pi,\pi],
  \label{eq:generaldomain} 
\end{equation}
where 
\begin{equation}
  r(\theta) = \sqrt{ \cos 2\theta + \sqrt{ 1.1 - \sin^{2}
      2\theta } }.
  \label{eq:peanut} 
\end{equation}
Recall, for each $y^{\star}$, the boundary and domain is rotated from the $(\theta, \phi)$ coordinate system to the $(s,t)$ coordinate {system} so that $y^{\star} = y(0,\cdot)$.
We present results for all three methods at one $y^{\star}$ with varying $\epsilon$ ($\ell = 1$) values, as shown in Fig. \ref{fig:other_quad_methods}A. 
 We consider the exact solution, from Beale {\it et al.}~\cite{beale2016simple},
 \begin{equation}
   u_{\text{ex}}(x_1,x_2,x_3) =  e^{x_{3}} ( \sin
   x_{1} + \sin x_{2} ),
   \label{eq:harmonicfctothermethods}
 \end{equation}
 with $(x_1,x_2,x_3)$ denoting an ordered triple in a Cartesian coordinate system.
 We substitute \eqref{eq:harmonicfctothermethods} evaluated {on} the boundary  and its normal derivative evaluated on the
 boundary into \eqref{eq:repres-int-sub} and
 compute numerical approximations of this harmonic function using \eqref{eq:genquadrature}.
\begin{figure}[htb]
\centering
\includegraphics[width=\linewidth]{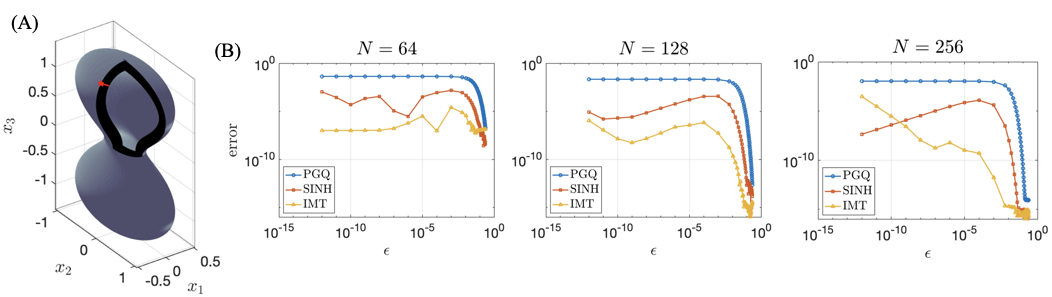}
\caption{ {(A) The peanut-shaped domain and the point $y^{\star} = (-0.4349, 0, 1.1819)$ (red star) with varying $\epsilon$ values (red line) at which we compute the error in evaluating \eqref{eq:repres-int-sub} using \eqref{eq:genquadrature} with three different methods. (B) The logarithmic error as a function of $\epsilon$ for the three quadrature rules used, PGQ, SINH, and IMT, to evaluate \eqref{eq:genquadrature} for three fixed values of $N$.}}
\label{fig:other_quad_methods}
\end{figure}}

{In Fig. \ref{fig:other_quad_methods}B, we present the error with respect to $\epsilon$  when evaluating \eqref{eq:repres-int-sub} with the three numerical methods  for different fixed $N$ values, $N =64$, $128$, and $256$. The PGQ is a method  often used for evaluating double- and single-layer potentials in three dimensions but does not treat the close evaluation problem directly. This is what we observe here as the PGQ has large errors for all values of fixed $N$ as  $\epsilon \to 0^{+}$. The other two quadrature methods, SINH and IMT, are specifically designed for nearly singular integrals and choose abscissas such that more quadrature points are located near the nearly singular point, $s = 0$, see Fig. \ref{fig:gauss-grid}. We refer the reader to \cite{iri1987certain, johnston2005sinh, robinson1981algorithm} for details on how these  abscissas and corresponding weights are chosen.  We point out the implementation of these two methods pose some challenges, including that the IMT quadrature requires the use of a table for the abscissas and weights and the SINH quadrature rule requires using $\epsilon$ to determine appropriate abscissas and weights.
We first observe in Fig. \ref{fig:other_quad_methods} that the SINH and IMT methods do improve the error in the close evaluation problem and the errors are not monotonic in $\epsilon$. The SINH quadrature for a fixed $N$ that is large, $N = 256$, does exhibit some convergence but is at a rate significantly less than $O(\epsilon)$. }

{ Furthermore, for these methods we do not see convergence for a fixed $\epsilon$ as $N$ increases from $64$ to $256$. For these moderate $N$ values, the error is dominated by the close evaluation error. Note all of these methods should be convergent as $N \to \infty$ for a fixed $\epsilon$ provided that the number of quadrature points, $N$, is sufficiently large. The challenge with these problems is that the value of a sufficiently large $N$ is necessarily huge for these nearly singular integrals. 
}

{In this paper, we develop a numerical method and quadrature rule that takes into account the underlying behavior
of these integrals which we investigate through a local analysis. We develop a method that does well for a fixed N, in
particular that is not large, as $\epsilon \to 0^{+}$. Furthermore, the implementation of this method is straightforward, similar to PGQ.}

\section{Local analysis}
\label{sec:local-analysis}

Since the sharply peaked behavior of \eqref{eq:repres-int-sub}
about $y = y^{\star}$ is the major cause for error in its evaluation,
we analyze the contribution made by a small, fixed region about
$y = y^{\star}$ in the asymptotic limit, $\epsilon \to 0^{+}$. This
setting represents the situation in which the distance of the
evaluation point from the boundary,  $\epsilon \ell$, is
smaller than the discretization used by a numerical
integration method at a fixed resolution{, $N$, }for \eqref{eq:repres-int-sub}.  Through the following
local analysis, we obtain valuable insight into the asymptotic
behavior of \eqref{eq:repres-int-sub} as $\epsilon \to 0^{+}$.

Let $B_{\delta}$ denote a fixed, small portion of $B$ that contains
the point $y^{\star}$, {where $(s,t) \in [0, \delta] \times [-\pi, \pi]$}, and let
\begin{eqnarray}
  U^{\delta}(y^{\star};\epsilon) &=& \frac{1}{4 \pi} \int_{B_{\delta}}
  \frac{n(y) \cdot ( y^{\star} - y - \epsilon \ell n^{\star})}{|
    y^{\star} - y - \epsilon \ell n^{\star} |^{3}} \left[ \mu(y) -
    \mu(y^{\star}) \right] 
  \mathrm{d}\sigma_{y}   \label{eq:Ulocal} \\ 
&=& { \frac{1}{4 \pi} \int_{-\pi}^{\pi}
  \int_{0}^{\delta}
  \frac{\tilde{n}(s,t) \cdot ( y(0,\cdot) - y(s,t) - \epsilon \ell
   n^\star)}{| y(0,\cdot) - y(s,t) - \epsilon \ell
      n^\star|^{3}}
  J(s,t) \left[ \tilde{\mu}(s,t) - \tilde{\mu}(0,\cdot) \right]
  \sin(s) \,  \mathrm{d}s \mathrm{d}t,}
  \label{eq:Udelta}
\end{eqnarray}
denote the contribution to the double-layer potential made by
$B_{\delta}$, and let
\begin{eqnarray}
  V^{\delta}(y^{\star};\epsilon) &=& \frac{1}{4\pi} \int_{B_{\delta}}
  \frac{1}{|y^{\star} - y - \epsilon \ell n^{\star}|} \rho(y) \,
  \mathrm{d}\sigma_{y}   \label{eq:Vlocal} \\ 
  &=& { \frac{1}{4\pi} \int_{-\pi}^{\pi}
  \int_{0}^{\delta} \frac{J(s,t)}{|y(0,\cdot) - y(s,t) - \epsilon \ell
     n^\star|} \tilde{\rho}(s,t)   \sin(s) \, \mathrm{d}s
  \mathrm{d}t,}
  \label{eq:Vdelta}
\end{eqnarray}
denote the contribution to the single-layer potential made by
$B_{\delta}$. { Here, the spherical Jacobian, $\sin(s)$, is explicitly included in the
integral and $J(s,t) = | \partial_{s} y(s,t) \times \partial_{t} y(s,t)
|/\sin(s) =  \tilde{J}(s,t)/\sin(s)$.   The global Jacobian, $\tilde{J}(s,t) = J(s,t) \sin(s)  =  | \partial_{s} y(s,t) \times \partial_{t} y(s,t)|$, remains bounded.}
{In the results that follow,  where \eqref{eq:Udelta} and \eqref{eq:Vdelta} are analyzed in
the limit as $\epsilon \to 0^{+}$, we provide insight into the close
evaluation of the double- and single-layer potentials which leads to a numerical method.
Layer potentials in spherical coordinates are well known \cite{folland1995introduction}, but their asymptotic behavior, especially using the expressions \eqref{eq:Udelta}-\eqref{eq:Vdelta}, is not well studied.}
\subsection{Local analysis of the double-layer potential}

The kernel in \eqref{eq:Udelta} is given by
\begin{equation}
  K(s,t;\epsilon) = \frac{\tilde{n}(s,t) \cdot ( y_{d}(s,t) -
    \epsilon \ell n^{\star} ) }{| y_{d}(s,t) - \epsilon \ell
    n^{\star} |^{3}} J(s,t),
  \label{eq:K-kernel}
\end{equation} 
with $y_{d}(s,t) = y(0,\cdot) - y(s,t)$. To study the asymptotic
behavior of $K$ as $\epsilon \to 0^{+}$, we introduce the stretched
coordinate, $s = \epsilon S$. Note that
$\tilde{n}(\epsilon S, t) = n^{\star} + O(\epsilon)$ and
$y_{d}(\epsilon S, t) = -\epsilon S y_{s}(0,\cdot) + O(\epsilon^{2})$
where $y_{s}(0,\cdot) \neq 0$ is a vector that lies on the tangent plane
orthogonal to {$n^\star$}. {Recall, that this vector, $y_s(0,\cdot)$, is defined by averaging the $t$ coordinate using the spherical mean.} We find by expanding
$K(\epsilon S,t;\epsilon)$ about $\epsilon = 0$ that
\begin{equation}
  K(\epsilon S,t;\epsilon) = - 
  \frac{\ell J(0,\cdot)}{\epsilon^{2} ( S^{2} | y_{s}(0,\cdot) |^{2} +
    \ell^{2} )^{3/2}} + O(\epsilon^{-1}).
  \label{eq:kernel-behavior}
\end{equation}
There are two key observations about this leading behavior of
$K$. In the limit as $\epsilon \to 0^{+}$, $K = O(\epsilon^{-2})$
which characterizes the nearly singular behavior of the double-layer
potential. In addition, {the} leading behavior of $K$ given in
\eqref{eq:kernel-behavior} is independent of the azimuthal angle,
$t$. In other words, $K$ is azimuthally invariant about $y^{\star}$ in the
limit as $\epsilon \to 0^{+}$ to leading order.

When we substitute $s = \epsilon S$ into \eqref{eq:Udelta}, replace
$K$ by its leading behavior given in \eqref{eq:kernel-behavior}, and use  
 $\sin( \epsilon S) = \epsilon S + O(\epsilon^3)$, we find
\begin{equation}
  U^{\delta}(y^{\star};\epsilon) = -\frac{\ell J(0,\cdot)}{2}
  \int_{0}^{\delta/\epsilon} \left( \frac{S}{(S^{2} | y_{s}(0,\cdot) |^{2} +
    \ell^{2} )^{3/2}} + O(\epsilon) \right) \left[ \frac{1}{2\pi} \int_{-\pi}^{\pi} \left(
      \tilde{\mu}(\epsilon S,t)
      - \tilde{\mu}(0,\cdot)  \right) \mathrm{d}t \right] \mathrm{d}S \, .
  \label{eq:Udelta-intermediate}
\end{equation}
Note that the integration in $t$ in the square brackets above is the
average of $\tilde{\mu}(\epsilon S,t) - \tilde{\mu}(0,\cdot)$ over a
circle about $y^{\star}$.  Naively, it appears that because
$\tilde{\mu}(\epsilon S,t) - \tilde{\mu}(0,\cdot) = O(\epsilon)$ that
the integral in $t$ will be $O(\epsilon)$. However, we find that
\begin{align}
  \frac{1}{2\pi} \int_{-\pi}^{\pi} \left[ \tilde{\mu}(\epsilon S,t) -
  \tilde{\mu}(0,\cdot) \right] \mathrm{d}t
  &= \frac{1}{2\pi} \int_{0}^{\pi} \left[ \tilde{\mu}(\epsilon S,t)
    + \tilde{\mu}(\epsilon S,t - \pi) - 2 \tilde{\mu}(0,\cdot) \right]
    \mathrm{d}t \nonumber\\
  &= \frac{1}{2\pi} \int_{0}^{\pi} \left[ \tilde{\mu}(\epsilon S,t)
    + \tilde{\mu}(-\epsilon S,t) - 2 \tilde{\mu}(0,\cdot) \right]
    \mathrm{d}t \nonumber\\
  & = \epsilon^{2} S^{2} \frac{1}{2\pi} \int_{0}^{\pi} \partial_{s}^{2}
    \tilde{\mu}(0,\cdot) \, \mathrm{d}t  + O(\epsilon^4) \, .
    \label{eq:spherical-laplacian}
\end{align}
In fact, the averaging operation yields a result
that is $O(\epsilon^{2})$.  Here, we have used the regularity of
$\tilde{\mu}$ over the north pole to substitute
$\tilde{\mu}(\epsilon S,t - \pi) = \tilde{\mu}(-\epsilon
S,t)$. Furthermore, it has been shown by the
authors~\cite{ckk2018asymp} that
\begin{equation} \label{eq:laplacian_S}
  \frac{1}{2 \pi} \int_{0}^{\pi} \partial_{s}^{2}
    \tilde{\mu}(0,\cdot) \mathrm{d}t = \frac{1}{4} \varDelta_S
    \tilde{\mu}(0,\cdot),
\end{equation}
with $\varDelta_S \tilde{\mu}(0,\cdot)$ denoting the spherical Laplacian
of $\tilde{\mu}$ evaluated at the north pole. Substituting
\eqref{eq:spherical-laplacian} and \eqref{eq:laplacian_S}  into \eqref{eq:Udelta-intermediate}, we
obtain
\begin{equation}\label{eq:second_order_mu}
  U^{\delta}(y^{\star};\epsilon) = -\frac{\ell J(0,\cdot)}{8}
  \varDelta_S \tilde{\mu}(0,\cdot) \int_{0}^{\delta/\epsilon}
  \frac{\epsilon^{2} S^{3}}{(S^{2} | y_{s}(0,\cdot) |^{2} + \ell^{2}
    )^{3/2}} \, \mathrm{d}S + O(\epsilon^3)\, .
\end{equation}
In the asymptotic limit corresponding to
$0 < \epsilon \ll \delta \ll 1$, we find that
\begin{equation}
  \int_{0}^{\delta/\epsilon}
  \frac{\epsilon^{2} S^{3}}{(S^{2} |
    y_{s}(0,\cdot) |^{2} + \ell^{2} )^{3/2}} \, \mathrm{d}S  
=\delta \frac{\epsilon}{|y_{s}(0,\cdot)|^{3}} + O(\epsilon^2).
  \label{eq:script-I}
\end{equation}
It follows that
\begin{equation}
  U^{\delta}(y^{\star};\epsilon,\delta) = - \delta \frac{\epsilon \ell
    J(0,\cdot)}{8 |y_{s}(0,\cdot)|^{3}} \varDelta_S
  \tilde{\mu}(0,\cdot) + O(\epsilon^2).
  \label{eq:DLP-local}
\end{equation}
This result gives the leading behavior for $U^{\delta}$. The results
\eqref{eq:kernel-behavior}, \eqref{eq:script-I} and \eqref{eq:DLP-local} provide the following
insights about the close evaluation problem.
\begin{itemize}

\item The kernel is azimuthally invariant about $y^{\star}$ as
  $\epsilon \to 0^{+}$ to leading order.

\item Integrating the difference $\mu(y) - \mu(y^{\star})$ over a
  closed circuit surrounding $y^{\star}$ introduces an averaging
  operation that effectively smooths that difference as shown in
  \eqref{eq:spherical-laplacian}.
  
\item It follows that the leading order behavior gives
  $U^{\delta}(y^{\star};\epsilon) = O(\epsilon)$ as
  $\epsilon \to 0^{+}$.
 
\end{itemize}

\begin{figure}[htb]
  \centering
  \includegraphics[width=0.95\linewidth]{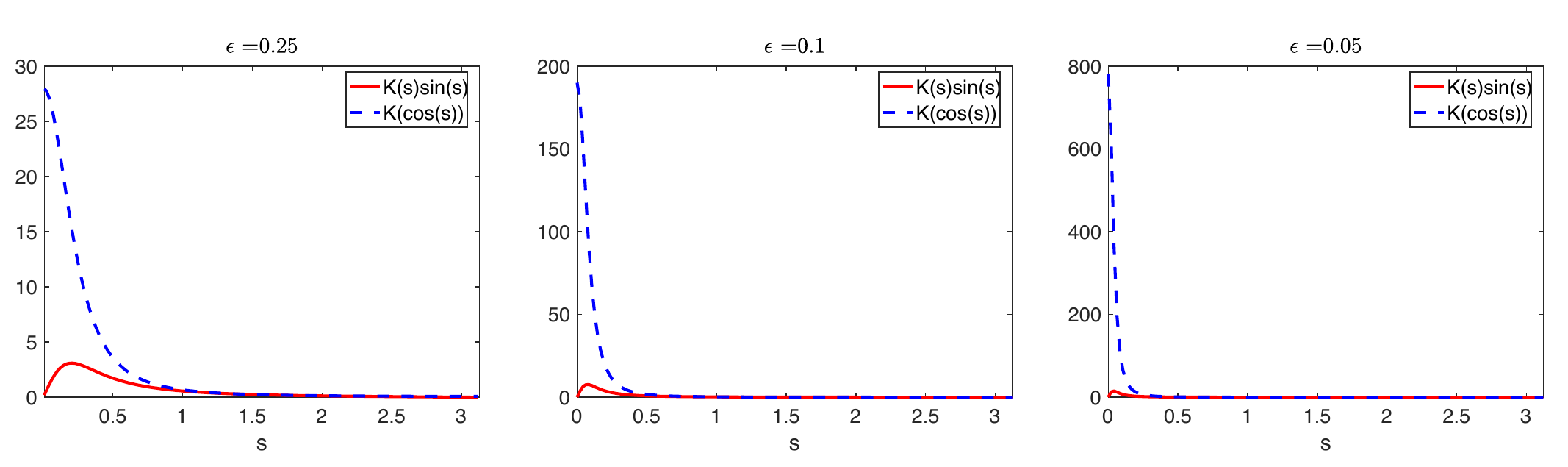}
  \caption{Plots of $K(s,t)\sin(s)$ (red {solid}) and
    $K(z,t) = K(\cos(s),t) $ (blue {dashed}) for the case in
    which the boundary is the unit sphere for different values of
    $\epsilon\ell $ ($\ell = 1$). For the unit sphere, we have
    $K(s,t)\sin(s) = K(s)\sin(s) = (2 \epsilon - \epsilon^{2}) \sin(s) \left[ 2 (1 -
      \epsilon) (1 - \cos s ) + \epsilon^{2} \right]^{-3/2}$ and
    $K(z,t) = K(z) = (2 \epsilon - \epsilon^{2}) \left[ 2 (1 - \epsilon) (1 -
      z ) + \epsilon^{2} \right]^{-3/2}$. 
     }
\label{img:compa}
\end{figure}

For integrals in spherical coordinates, it is common to substitute $s$ by
$z = \cos(s)$ 
and use
Gauss-Legendre quadrature to integrate in $z$
\cite{atkinson1982numerical}. {This is the PGQ method presented in Section \ref{sec:prior}}. In this $(z,t)$-spherical coordinate
system \eqref{eq:Udelta} becomes
\begin{align}
  \tilde{U}^{h}(y^{\star};\epsilon)
  &= \frac{1}{4 \pi} \int_{-\pi}^{\pi} \int_{1-h}^{1}
    \frac{\hat{n}(z,t) \cdot ( y(1,\cdot) - y(z,t) - \epsilon
    \ell n^\star)}{| y(1,\cdot) - y(1,t) - \epsilon \ell n^\star|^{3}}
    \hat{J}(z,t) \left[ \hat{\mu}(z,t) -
    \hat{\mu}(1,\cdot) \right] \mathrm{d}z \mathrm{d}t
    \nonumber\\
  &= \frac{1}{4 \pi} \int_{-\pi}^{\pi} \int_{1-h}^{1} K(z,t) \left[
    \hat{\mu}(z,t) - \hat{\mu}(1,\cdot) \right]
    \mathrm{d}z \mathrm{d}t,
  \label{eq:Udelta_tilde}
\end{align}
with $\hat{n}(z,t) = n(y(z,t))$, $\hat{\mu}(z,t) = \mu(y(z,t))$,
$h = 1 - \cos(\delta)$, and
$\hat{J}(z,t) = | \partial_{z} y(z,t) \times \partial_{t} y(z,t)|$. To study
the asymptotic behavior of $K(z,t)$, we substitute
$z = 1- \epsilon Z$, expand $K(1 - \epsilon Z, t ; \epsilon)$ about
$\epsilon = 0$ and find that
\begin{equation}
  K(1- \epsilon Z,t;\epsilon) = \displaystyle - \frac{\ell
  \hat{J}(1,\cdot)}{\epsilon^{2} ( Z^{2} | y_{z}(1,\cdot) |^{2}
  + \ell^{2} )^{3/2}} + O(\epsilon^{-1}).
\end{equation}
This asymptotic behavior is the same as 
 in
\eqref{eq:kernel-behavior}. However, integrating in the
$(s,t)$-coordinate system includes the factor $\sin(s)$, which
effectively reduces the peak about $s = 0$ that causes the nearly
singular behavior of this integral. To show this, we plot in
Fig.~\ref{img:compa} a comparison between $K(s,t)\sin(s)$ and
$K(z,t) = K(\cos(s),t) $ for the 
 case in which the boundary
is the unit sphere. {Note,  the factor $\sin(s)$, used explicitly, contributes to remove the nearly peaked behavior in the integrand. Furthermore,} in the asymptotic analysis above, this factor of
$\sin(s)$ contributes an $O(\epsilon)$ factor leading to
\eqref{eq:DLP-local}. {This observation leads us to choose a different substitution to go from  $s \in [0, \pi]$ to $z \in [-1,1]$, $z = \frac{2}{\pi}s -1$, which allows us to naturally take advantage of the explicit $\sin(s)$ factor. Then, we can apply a quadrature method in $z \in [-1,1]$ that addresses the close evaluation problem.}

\subsection{Local analysis of the single-layer potential}

\begin{figure}[h!]
  \centering
  \includegraphics[width=0.95\linewidth]{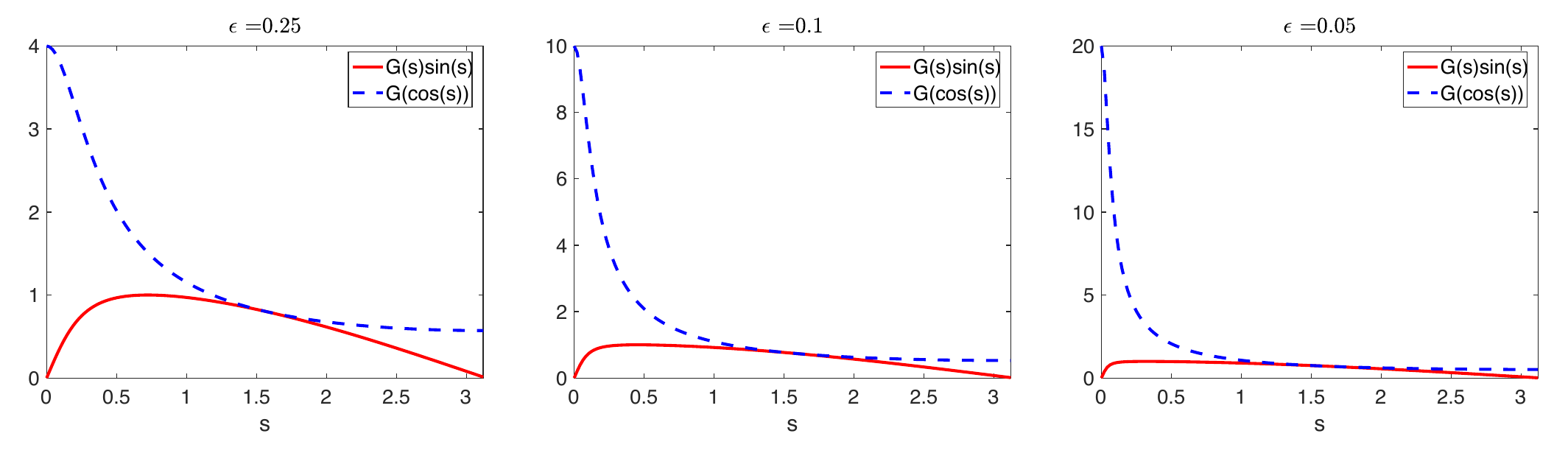}
  \caption{Plots of $G(s,t)\sin(s)$ (red {solid}) and
    $G(z,t) = G(\cos(s),t) $ (blue {dashed}) for the case in
    which the boundary is the unit sphere for different values of
    $\epsilon\ell$ ($\ell = 1$).  For the unit sphere, we have
    $G(s,t)\sin(s) = G(s)\sin(s) = \sin(s) \left[ 2 (1 - \epsilon) (1 - \cos s ) +
      \epsilon^{2} \right]^{-1/2}$ and
    $G(z,t) = G(z) = \left[ 2 (1 - \epsilon) (1 - z ) + \epsilon^{2}
    \right]^{-1/2}$. 
 }
\label{img:compa2}
\end{figure}

We now study the asymptotic behavior of
$V^{\delta}(y^{\star};\epsilon)$ as $\epsilon \to 0^{+}$. The steps
that we take here follow those used above for
$U^{\delta}(y^{\star};\epsilon)$. In particular, we introduce the
stretched coordinate $s = \epsilon S$ and find that the kernel in
\eqref{eq:Vdelta}, which we denote by $G$, has the leading behavior,
\begin{equation}
  G(\epsilon S,t;\epsilon) = \frac{J(\epsilon S,t)}{| y_{d}(\epsilon
    S, t) - \epsilon \ell \tilde{n}(0,\cdot) |}=
  \frac{J(0,\cdot)}{\epsilon ( S^{2} | y_{s}(0,\cdot) |^{2} + \ell^{2}
    )^{1/2}}  + O(1).
  \label{eq:Gkernel-behavior}
\end{equation}
In the limit as $\epsilon \to 0^{+}$, we find that
$G = O(\epsilon^{-1})$ which characterizes its nearly singular
behavior.  Just as with the double-layer potential, we find that the
leading behavior of this kernel is independent of $t$ and hence is
azimuthally invariant about $y^{\star}$. When we substitute $s = \epsilon S$
into \eqref{eq:Vdelta} and replace $G$ by its leading behavior given
in \eqref{eq:Gkernel-behavior}, we find after integrating with respect
to $t$,
\begin{equation}
  V^{\delta}(y^{\star};\epsilon) = \frac{1}{2}
  J(0,\cdot) \tilde{\rho}(0,\cdot) 
  \int_{0}^{\delta/\epsilon} \frac{\epsilon S}{(S^{2} | y_{s}(0,\cdot) |^{2} +
    \ell^{2})^{1/2}} \, \mathrm{d}S + O(\epsilon^2) \, .
\end{equation}
In the asymptotic limit corresponding to
$0 < \epsilon \ll \delta \ll 1$, we find that
\begin{equation}
  \int_{0}^{\delta/\epsilon} \frac{\epsilon S}{( S^{2} |
    y_{s}(0,\cdot) |^{2} + \ell^{2} )^{1/2}} \, \mathrm{d}S =  \delta
  \frac{1}{|y_{s}(0,\cdot)|}  + O(\epsilon) \, .
\end{equation}
It follows that
\begin{equation}
  V^{\delta}(y^{\star};\epsilon) =\delta \frac{1}{2} \frac{J(0,\cdot)
    \tilde{\rho}(0,\cdot)}{|y_{s}(0,\cdot)|} + O(\epsilon).
  \label{eq:SLP-local}
\end{equation}
In contrast to $U^{\delta}$, we find that
the leading order behavior gives $V^{\delta}(y^{\star};\epsilon) = O(1)$ as $\epsilon \to 0^{+}$.\\
For the same reasons as discussed above for the double-layer potential, one finds that $G(s,t)\sin(s)$ is a smoother
integrand than $G(z,t) = G(\cos(s),t)$.  In Fig. \ref{img:compa2} we present
an example for the 
 case of the unit sphere. {Once again, this observation leads us to choose the substitution, $z = \frac{2}{\pi}s -1$.}

\subsection{Results from this local analysis}
\label{sec:localanlysisresults}

The contribution by $B_{\delta}$ to \eqref{eq:repres-int-sub} is given
by
\begin{equation}
  u^{\delta}(y^{\star} - \epsilon \ell n^{\star}) = u(y^{\star}) -
  U^{\delta}(y^{\star};\epsilon) + V^{\delta}(y^{\star};\epsilon).
\end{equation}
Using the results for the leading behaviors of $U^{\delta}$ and
$V^{\delta}$ from above, we expect the error made by the leading
behavior of $U^{\delta}$ given by \eqref{eq:DLP-local} to be
$O(\epsilon^{2})$ and the error made by the leading behavior of
$V^{\delta}$ given by \eqref{eq:SLP-local} to be
$O(\epsilon)$. Consequently, the error made by the leading behavior of
$u^{\delta}$ is $O(\epsilon)$.

The asymptotic behavior computed above provides valuable insight into
 {developing an accurate numerical method for}  \eqref{eq:repres-int-sub}. For any
numerical integration method used to compute \eqref{eq:repres-int-sub}
at close evaluation points, large errors are likely to occur when the
normal distance from $B$, given by $\epsilon \ell$, is smaller than the
fixed discretization length, 
{$\pi/N$, see PGQ in Fig. \ref{fig:other_quad_methods}}. {In} that case, the kernels 
become sharply peaked about $y^{\star}$ which
cannot be adequately resolved on the fixed boundary mesh.

The local analysis above indicates that it is advantageous in this
setting to consider: {(i) a rotated spherical} coordinate system defined with respect
to $y^{\star}$ {that enhances the
asymptotic behavior of the double- and single-layer potentials, and (ii) the use of the spherical Jacobian on the unit sphere, $\sin (s)$, explicitly that  yields a smoother integrand}. 
{Following the above guidelines, one can guarantee an approximation of \eqref{eq:repres-int-sub} that converges linearly with respect to distance from the boundary, $\epsilon$, for a fixed resolution, $N$.}
\section{Numerical method for the close evaluation of layer potentials}
\label{sec:numerics}

We present a new numerical method for computing \eqref{eq:repres-int-sub}
that makes use of the results from the the asymptotic analysis in
Section \ref{sec:local-analysis}. 
By 
{rotating the integrals as described above}, azimuthal integration in this rotated
coordinate system will naturally achieve {the averaging operation that reveals the asymptotic behavior of the layer potentials.} {Recall,} to
rotate the
 integrals in this way, we apply the 
  {tranformation matrix} 
described in \ref{sec:rotation}. The result of this rotation leads to 
\begin{equation}
  I(y^{\star}) = \frac{1}{4\pi} \int_{-\pi}^{\pi}
  \int_{0}^{\pi} {F}(s,t) \sin (s) \,\mathrm{d}s \mathrm{d}t,
  \label{eq:rotated-I3D}
\end{equation}
with
\begin{multline}
{
F(s,t ) = -
   \frac{\tilde{n}(s,t) \cdot ( y(0,\cdot) - y(s,t) - \epsilon \ell
   n^\star)}{| y(0,\cdot) - y(s,t) - \epsilon \ell
      n^\star|^{3}}
 J(s,t) \left[ \tilde{\mu}(s,t) - \tilde{\mu}(0,\cdot) \right]
  + \frac{1}{|y(0,\cdot) - y(s,t) - \epsilon \ell
     n^\star|} J(s,t) \tilde{\rho}(s,t)} \, . 
  \label{eq:integrand}
\end{multline}

Since the local analysis also benefited from the explicit
consideration of the factor of $\sin(s)$ appearing in the
integral, compare \eqref{eq:integrals} with {\eqref{eq:rotated-I3D}}, we choose to integrate with respect to the polar angle $s$
rather than the cosine of the polar angle as is typically done.  Let
$z_{i} \in (-1,1)$ for $i = 1, \cdots, N$ denote the $N$-point
Gauss-Legendre quadrature rule abscissas with corresponding weights,
$w_{i}$ for $i = 1, \cdots, N$. Let $t_{j} = -\pi + \pi (j - 1)/N$ for
$j = 1, \cdots, 2N$ denote an equi-spaced grid for $t \in [-\pi,\pi]$
for the periodic trapezoid rule. To compute \eqref{eq:rotated-I3D}
numerically, we introduce the 
{substitution} $s_{i} = \pi ( z_{i} + 1 ) / 2$
for $i = 1, \cdots, N$ and use the following quadrature rule,
\begin{equation}
  I(y^{\star}) \approx I^{N}(y^{\star}) = \frac{\pi}{8 N} \sum_{i =
    1}^{N} \sum_{j = 1}^{2N} w_{i} \sin( s_{i} )
F(s_{i},t_{j}).
  \label{eq:quadrature}
\end{equation}
In \eqref{eq:quadrature} a factor of $\pi/2$ appears due to the
scaling of the quadrature weights, $w_{i}$ for $i = 1, \cdots, N$,
that is required for the 
{substitution} from $z$ to $s$. 

\begin{figure}[htb]
  \centering
\includegraphics[width=0.5\linewidth]{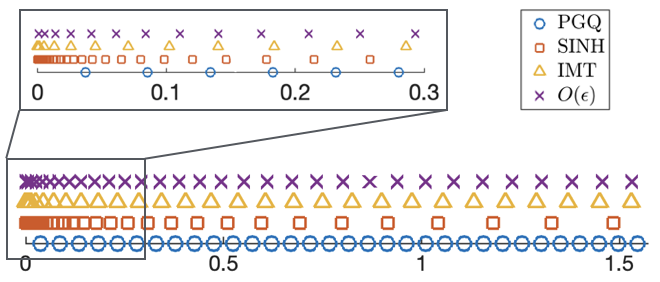}
  \caption{{   Comparison of the abscissas for PGQ, SINH quadrature, IMT quadrature, and our new method, denoted $O(\epsilon)$. 
     Each method has a resolution of $N =64$ between $0$ and $\pi$. Here, we present the quadrature points from $[0, \pi/2]$. and show a close up from $[0, 0.3]$. The abscissas location for the SINH quadrature is dependent on $\epsilon$ and we choose $\epsilon = 10^{-7}$ to show here.}
   }
  \label{fig:gauss-grid}
\end{figure}

The quadrature rule given in \eqref{eq:quadrature} is a modification
of the {PGQ rule by
Atkinson~\cite{atkinson1982numerical} presented in Section \ref{sec:prior}}. 
There are two factors that make
integrating with respect to $s$ more effective than integrating with
respect to $\cos(s)${: (1)} 
the inclusion of the factor $\sin(s)$
effectively smooths the peaks of the kernels in the double- and
single-layer potentials about $s = 0$ 
{and (2)} mapping the
Gauss-Legendre points from $(-1,1)$ to $(0,\pi)$ clusters the
quadrature points about $s = 0$ where the kernels for the double- and
single-layer potentials are peaked.
{In
Fig.~\ref{fig:gauss-grid}, we compare $s_i = \pi (z_{i} + 1) /2$, $i = 1, \dots, 64$, denoted $O(\epsilon)$, and
$\cos^{-1}(z_{i})$, the PGQ method presented in Section~\ref{sec:prior},} where {$z_{i}$ are} the $64$-point Gauss-Legendre
quadrature abscissas. {The two other quadrature methods presented in Section~\ref{sec:prior}, SINH and IMT, are designed to effectively 
cluster points near $s = 0$ to resolve the nearly singular behavior. We also include the SINH and IMT abscissas in Fig.~\ref{fig:gauss-grid}.
The new method and the IMT quadrature cluster in a similar way, while the SINH quadrature points cluster more, and the PGQ points do not cluster, as expected. 
Note that since the SINH abscissas depend on $\epsilon$, one moderate choice, $\epsilon = 10^{-7}$, is presented here. 
To develop a method that successfully clusters to resolve these integrals, the nearly peaked region should have more quadrature points, yet not too many points should be removed 
from other regions.  Observe that the SINH method is decreasing resolution quite significantly away from the peaked region and this leads to larger errors than IMT, as shown in Fig. \ref{fig:other_quad_methods}B.
We have chosen to use the Gauss-Legendre quadrature points in our method, based on PGQ, but one could also choose to use the IMT quadrature points with the mapping we suggest, $s_i = \pi (z_{i} + 1) /2$. We choose to modify the PGQ as that is a commonly used quadrature and does not require a tabulation as IMT quadrature does.}


{We also note that} since the Gauss-Legendre quadrature is an open rule, it does not
include the end points $z = \pm 1$. Consequently, the mapping
$s = \pi ( z + 1 ) / 2$ does not include the end points $s = 0$ and
$s = \pi$. {Therefore, the} quadrature rule
\eqref{eq:quadrature} does not require the explicit computation of
{$F(0,\cdot)$} which is the peak responsible for the nearly singular nature of this integral. 

The quadrature rule given in
\eqref{eq:quadrature} is effective for computing
\eqref{eq:repres-int-sub} because it excludes the need to evaluate the
function to be integrated at its peak, {smoothes the integrand (vanishing at the close evaluation point),} clusters the quadrature points near the peak{ed} behavior, and also includes the correct
averaging operation introduced by azimuthal integration in the rotated
coordinate system. 
{ Naively implementing this method to approximate the representation formula \eqref{eq:repres-int-sub}  for one evaluation point has the same computational cost as the PGQ method, $O(N^2)$. This is quite expensive, especially due to the rotation that is required for these methods. We have not focused in this paper on fast implementations but do believe that it is possible to speed up this method using ideas that have been previously developed including the fast multipole method \cite{FMM}, the spherical grid rotations of Gimbutas and Veerapaneni \cite{gimbutas2013fast}, and considering surface patches rather than the entire surface. }

\section{{Extension to $O(\epsilon^2)$}}
\label{sec:extension}
When quadrature rule \eqref{eq:quadrature} is adequately resolved for
$F$ defined in \eqref{eq:integrand}, we expect it to exhibit an
$O(\epsilon)$ error in the limit as $\epsilon \to 0^{+}${,} consistent
with the local analysis presented {in Section \ref{sec:local-analysis}}.  To obtain a smaller error,
we need to address the $O(\epsilon)$ error produced by the
single-layer potential. Instead of considering the single-layer potential directly, we
consider the asymptotic approximation{,}
\begin{equation}
  \frac{1}{4\pi} \int_{B}
  \frac{\rho(y)}{|y^{\star} - y - \epsilon \ell n^{\star}|} \,
  \mathrm{d}\sigma_{y} = \frac{1}{4\pi} \int_{B}
  \frac{\rho(y)}{|y^{\star} - y|} \,
  \mathrm{d}\sigma_{y} + \epsilon \ell \frac{1}{4\pi} \int_{B}
  \frac{n^{\star} \cdot ( y^{\star} - y )}{|y^{\star} - y|^{3}} \rho(y) \,
  \mathrm{d}\sigma_{y} - \frac{\epsilon \ell}{2} \rho(y^\star) + O(\epsilon^{2})\, ,
  \label{eq:SLP-perturbation}
\end{equation}
resulting from expanding the single-layer potential about
$\epsilon = 0$,  and using the boundary integral equation \cite{guenther1996partial}{,}
\[\frac{1}{4\pi} \frac{\partial}{\partial {n^\star}} \int_{B}
  \frac{\rho(y)}{|y^{\star} - y|} \,
  \mathrm{d}\sigma_{y} = - \frac{1}{2} \rho(y^\star) + \frac{1}{4\pi} \int_{B}
  \frac{n^{\star} \cdot ( y^{\star} - y )}{|y^{\star} - y|^{3}} \rho(y) \,
  \mathrm{d}\sigma_{y} .\] The asymptotic approximation given {in}
\eqref{eq:SLP-perturbation} is a sum of weakly singular integrals
defined on $B$. Unlike the double-layer potential, the
single-layer potential is continuous as $\epsilon \to 0^{+}$.  We
consider
 \eqref{eq:repres-int-sub} with the single-layer
potential replaced with \eqref{eq:SLP-perturbation},
\begin{multline}
  u(y^{\star} - \epsilon \ell n^{\star}) = u(y^{\star}) - \frac{1}{4
    \pi} \int_{B} \frac{n(y) \cdot ( y^{\star} - y - \epsilon \ell
    n^{\star})}{| y^{\star} - y - \epsilon \ell n^{\star} |^{3}}
  \left[ \mu(y) - \mu(y^{\star}) \right] \mathrm{d}\sigma_{y}\\
  + \frac{1}{4\pi} \int_{B} \frac{1}{|y^{\star} - y |} \rho (y) \,
  \mathrm{d}\sigma_{y} + \epsilon \ell \frac{1}{4\pi} \int_{B}
  \frac{n^\star \cdot ( y^{\star} - y )}{| y^{\star} - y |^{3}} \rho
  (y) \, \mathrm{d}\sigma_{y} - \frac{\epsilon \ell}{2} \rho(y^\star) {+ O(\epsilon^2)}.
  \label{eq:repres-int-sub-asympt}
\end{multline}
To compute a numerical approximation of the integrals in 
\eqref{eq:repres-int-sub-asympt}, we apply the same quadrature rule
described {in Section \ref{sec:numerics}} to the function
\begin{multline}
{
F(s,t ) = -
   \frac{\tilde{n}(s,t) \cdot ( y(0,\cdot) - y(s,t) - \epsilon \ell
   n^\star)}{| y(0,\cdot) - y(s,t) - \epsilon \ell
      n^\star|^{3}}
 J(s,t) \left[ \tilde{\mu}(s,t) - \tilde{\mu}(0,\cdot) \right] }\\
 {   + \frac{1}{|y(0,\cdot) - y(s,t) |} J(s,t) \tilde{\rho}(s,t) 
 + \epsilon \ell
  \frac{n^\star \cdot \left[
      y(0,\cdot) - y(s,t)
    \right]}{|y(0,\cdot) - y(s,t) |^3}
  J(s,t) \tilde{\rho}(s,t) \, .} 
  \label{eq:modified-integrand}
\end{multline}

Provided that quadrature rule \eqref{eq:quadrature} is adequately
resolved for \eqref{eq:modified-integrand}, we expect it to exhibit an
$O(\epsilon^{2})$ error as $\epsilon \to 0^{+}$.

\section{Numerical results}
\label{sec:examples}

 We present numerical results for the solution of Laplace's equation using representation formula \eqref{eq:repres_int} for close evaluation points. The new numerical method{,} detailed in
 Section \ref{sec:numerics}{,} is used to solve the formulas given in 
 \eqref{eq:repres-int-sub}  or \eqref{eq:repres-int-sub-asympt}. 
 We consider 
 exact solution {\eqref{eq:harmonicfctothermethods}, and proceed as in Section \ref{sec:prior}}.  

We test our method 
using two domains presented in Atkinson~\cite{atkinson1982laplace,
  atkinson1985algorithm}{: (1) the peanut-shaped domain, given by \eqref{eq:generaldomain} and \eqref{eq:peanut} and (2) a mushroom cap domain given by \eqref{eq:generaldomain} where }
%
%
 \begin{equation} 
  r(\theta) = 2 - \frac{1}{1 + 100 ( 1 - \cos\theta )^2}.
  \label{eq:mushroom} 
\end{equation} 
In the results that follow, we have set $\ell = 1$ for both of these
domains.

 \subsection{Results for peanut-shaped domain}


 \begin{figure}[h!]
  \centering
  \includegraphics[width=0.95\linewidth]{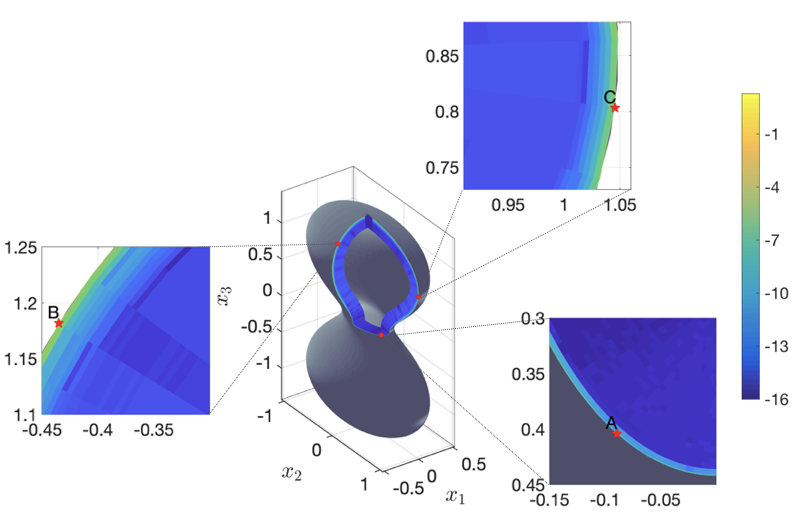}
  \caption{The logarithmic error for {the new numerical method with $N = 128$ applied to \eqref{eq:repres-int-sub}} when solving the close evaluation problem interior of the peanut-shaped domain. We have zoomed in around three points, A: $(-0.0894 , 0.4040, 0)$, B: $(-0.4349, 0, 1.1819)$, and C: $(0, 1.0456, 0.8032)$. Further data for these three points is shown in Figs.  \ref{img:peanut4} and \ref{img:peanut5}.}
\label{img:peanut2}
\end{figure}

\begin{figure}[h!]
 \centering
\includegraphics[width=0.32\linewidth]{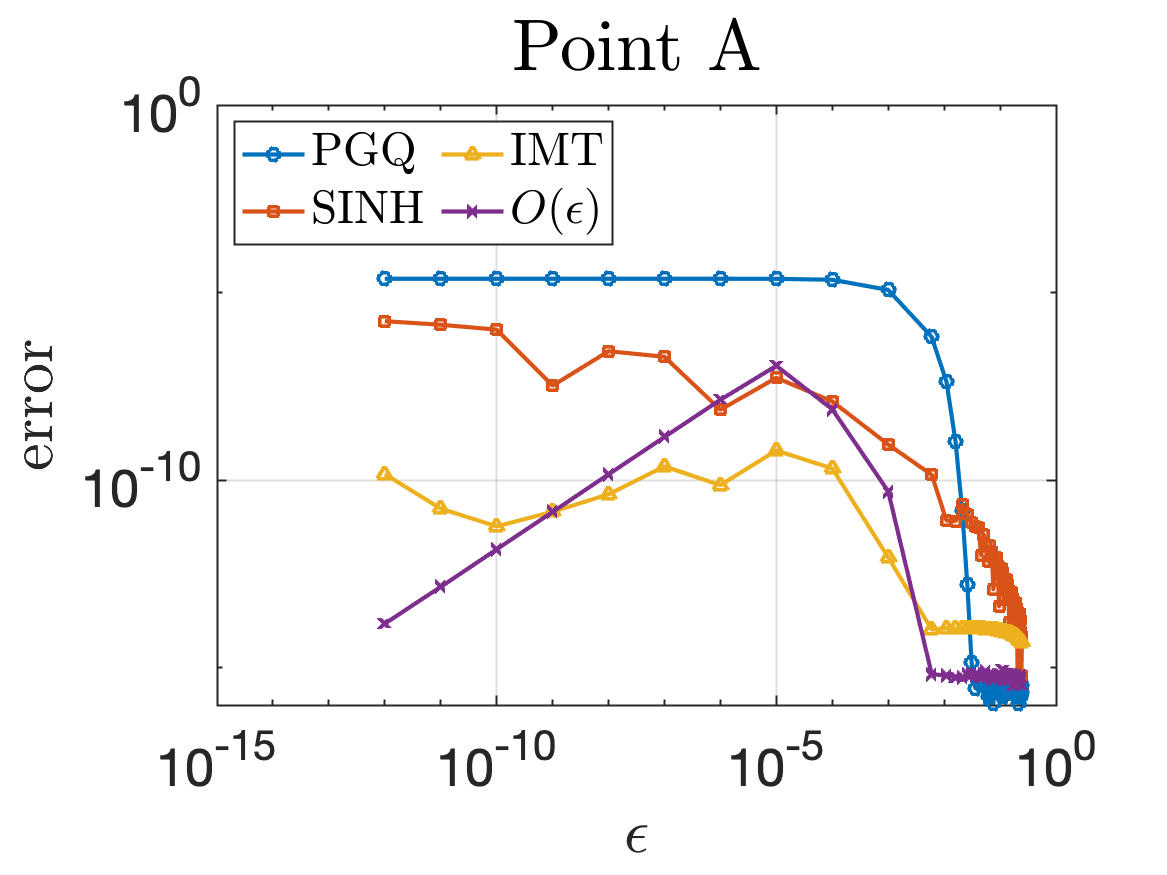}
\includegraphics[width=0.32\linewidth]{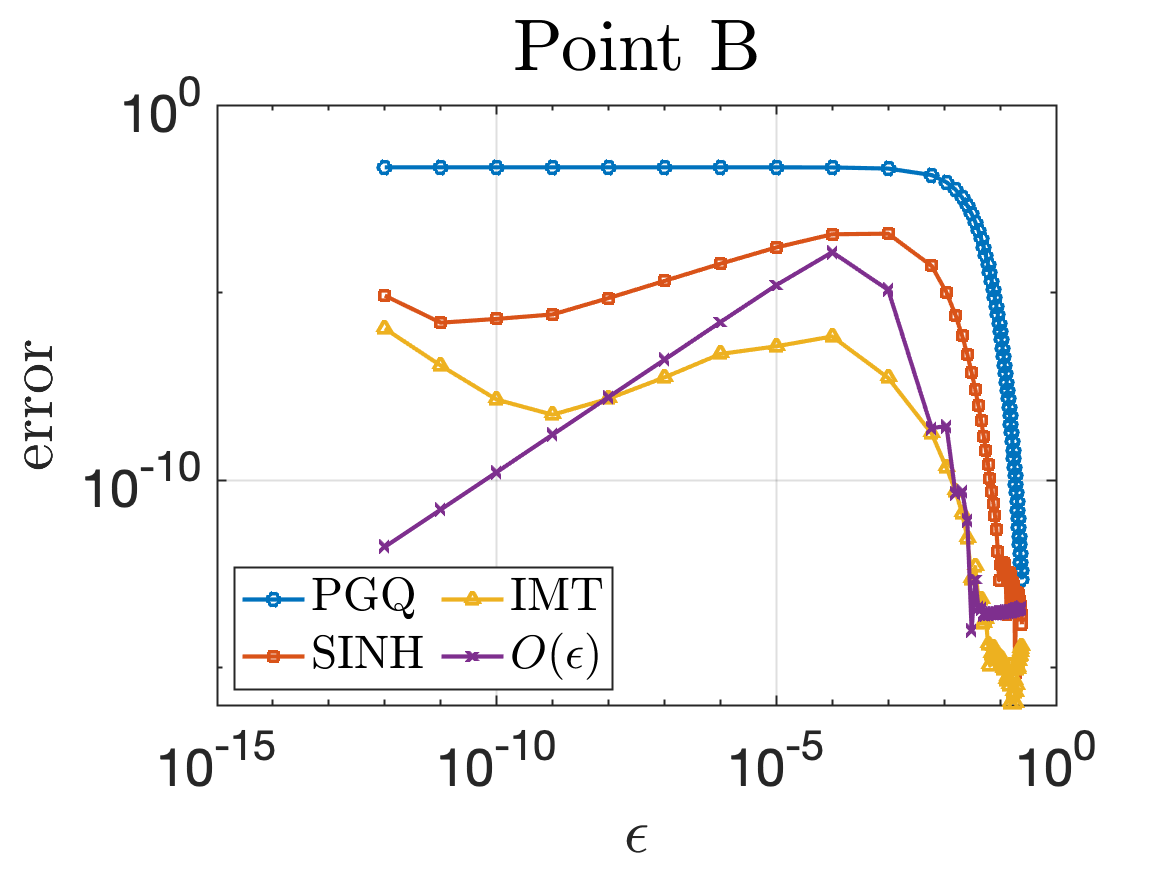}
\includegraphics[width=0.32\linewidth]{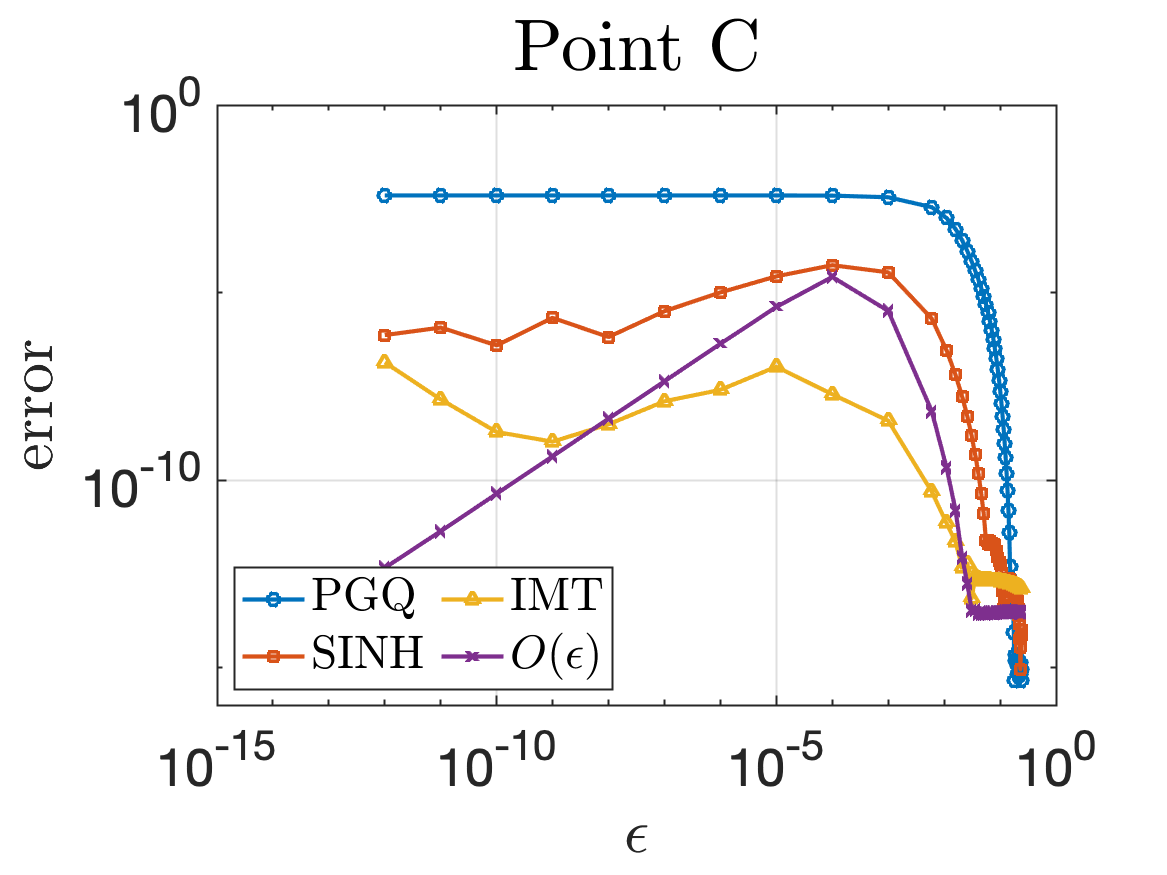}
  \caption{The logarithmic error for {the new numerical method applied to \eqref{eq:repres-int-sub} ($O(\epsilon)$), PGQ method, SINH method, and IMT method with $N = 128$}  at the three zoomed in points, A, B, and C of Fig. \ref{img:peanut2} when solving the close evaluation problem interior of the peanut-shaped domain as a function of $\epsilon$. }
\label{img:peanut4}
\end{figure}

\begin{figure}[h!]
 \centering
\includegraphics[width=0.32\linewidth]{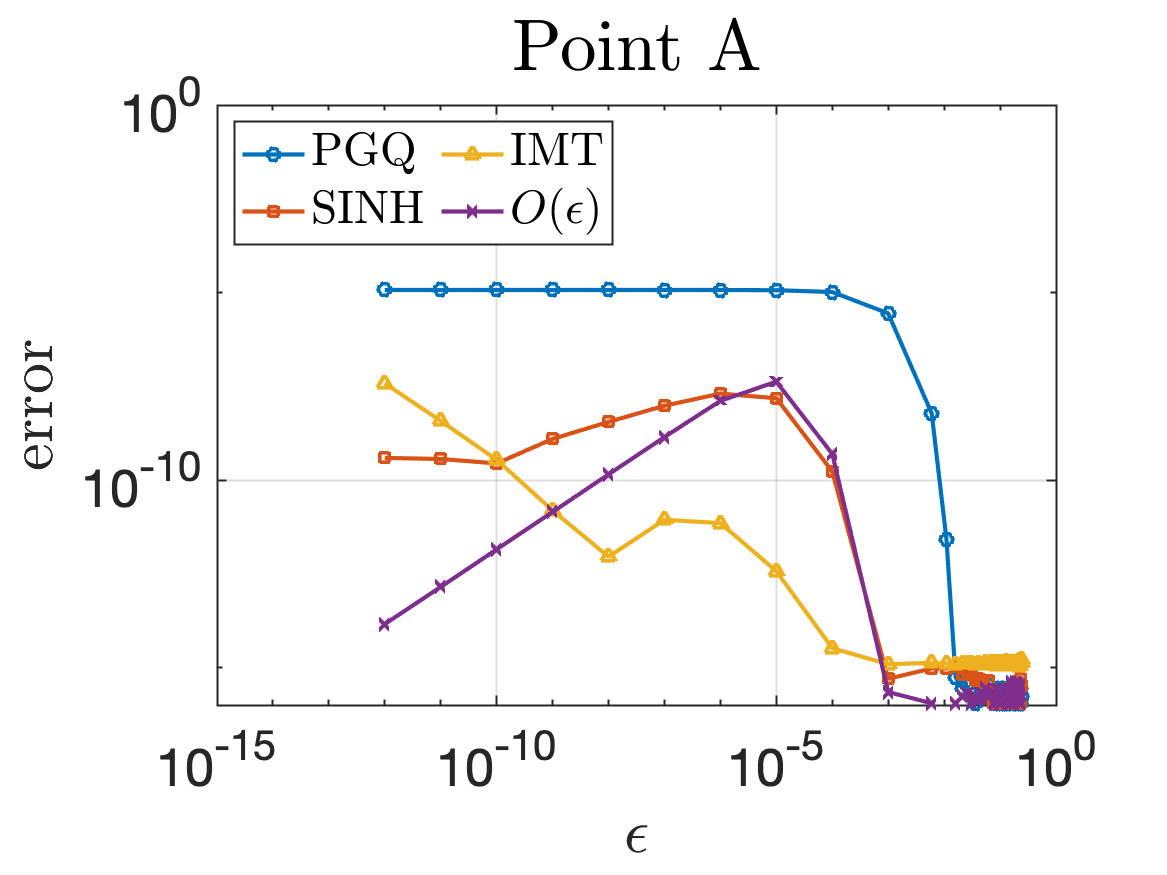}
\includegraphics[width=0.32\linewidth]{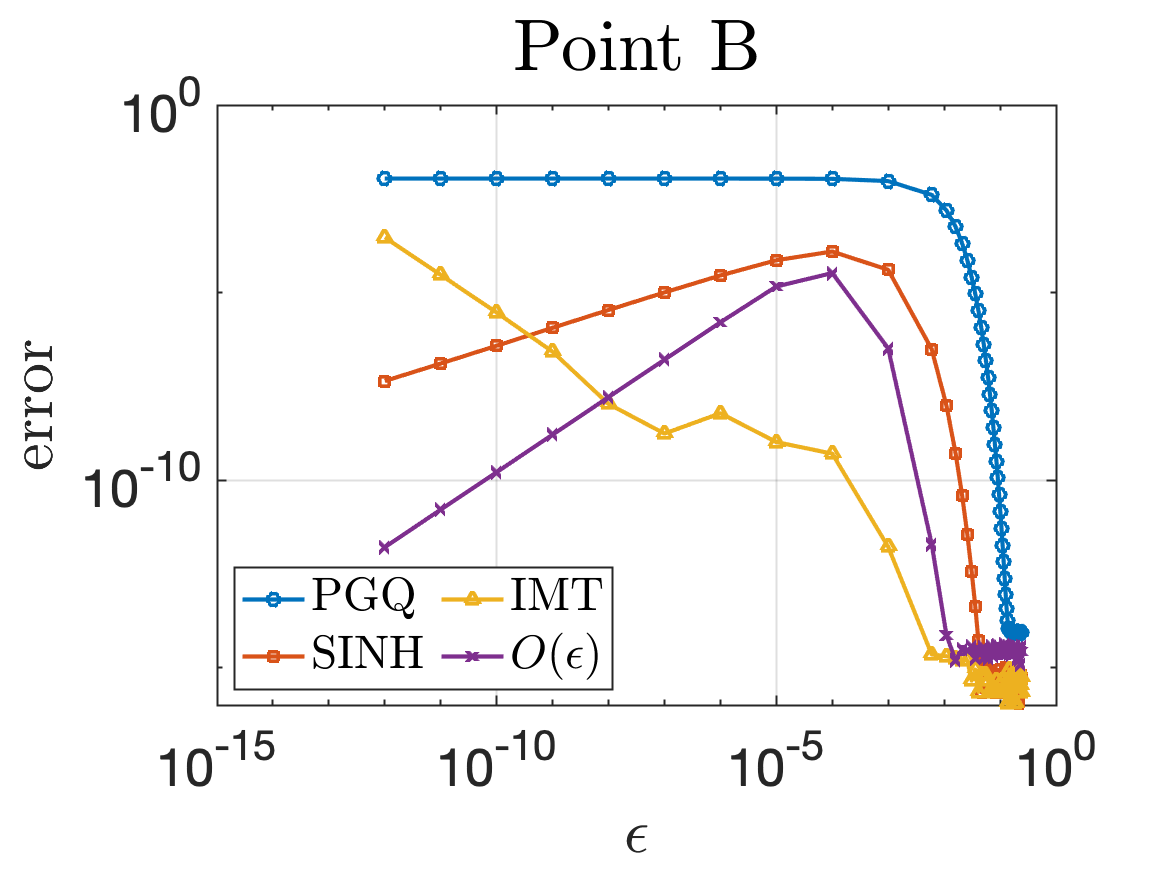}
\includegraphics[width=0.32\linewidth]{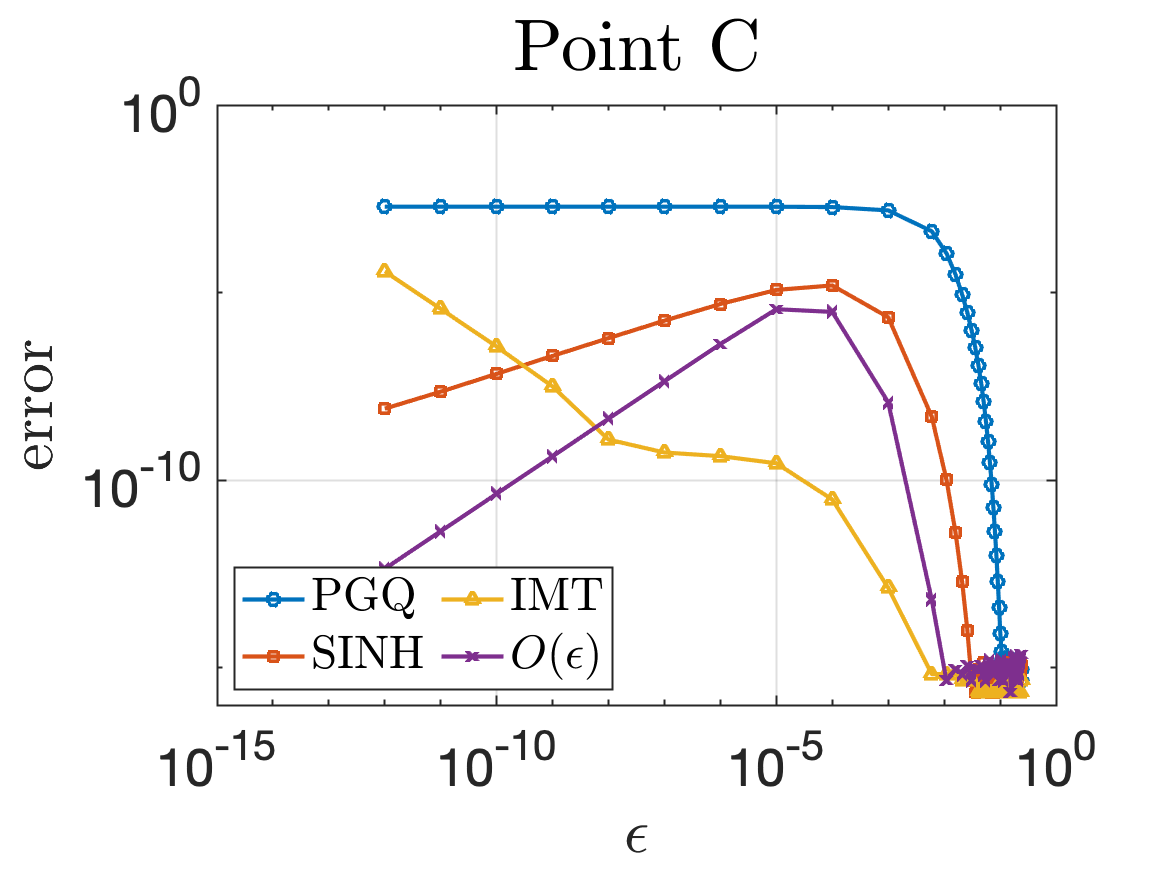}
  \caption{{The logarithmic error for the new numerical method applied to \eqref{eq:repres-int-sub} ($O(\epsilon)$), PGQ method, SINH method, and IMT method with $N = 256$ at the three zoomed in points, A, B, and C of Fig. \ref{img:peanut2} when solving the close evaluation problem interior of the peanut-shaped domain as a function of $\epsilon$. }}
\label{img:peanut5}
\end{figure}

In Fig. \ref{img:peanut2}, we present the error interior of the peanut-shaped domain when using our newly developed {numerical method with $N = 128$ to approximate \eqref{eq:repres-int-sub}}  
with the expectation of $O(\epsilon)$ error. In Fig. \ref{img:peanut4} {and Fig. \ref{img:peanut5}} we present the error as a function of $\epsilon$ for {this $O(\epsilon)$ method compared to the three prior methods that were introduced in Section \ref{sec:prior}, the} PGQ {method}, SINH {method}, and IMT {method}, 
{starting} at the three points (A, B, and C) labeled in Fig. \ref{img:peanut2}. {In Fig. \ref{img:peanut4}, all four methods use $N = 128$ and in Fig. \ref{img:peanut5}, all methods use $N = 256$.}  
 {We observe that  in all cases our method is $O(\epsilon)$ as $\epsilon \rightarrow 0^+$, as expected. Notably, this is even true for a moderate resolution of $N = 128$. Also, we note that in the case of higher resolution, $N = 256$, our method still has less error than all of the other three methods and further, the new method is the only one that has a consistent convergent behavior as $\epsilon \rightarrow 0^+$. }

 
\subsection{Results for mushroom cap domain}
 
\begin{figure}[h!]
  \centering
  \includegraphics[width=\linewidth]{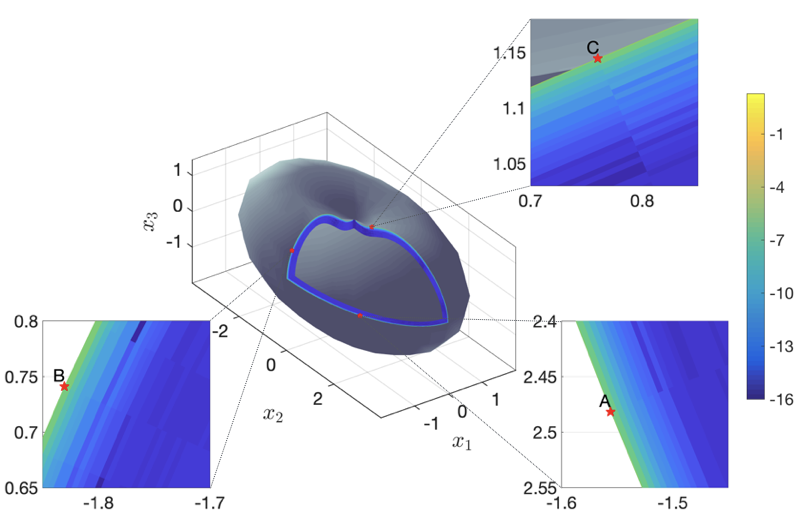}
  \caption{The logarithmic error for {the new numerical method with $N = 128$  applied to \eqref{eq:repres-int-sub}} when solving the close evaluation problem interior of the mushroom cap domain. We have zoomed in around three points, A: $(-1.5559, 2.4816, 0)$, B: $(-1.8307, 0, 0.7412)$, and C: $(0, 0.7601, 1.1446)$. Further data for these three points is shown in Fig.  \ref{img:mushroom4}.  }
\label{img:mushroom2}
\end{figure}

\begin{figure}[h!]
  \centering
\includegraphics[width=0.32\linewidth]{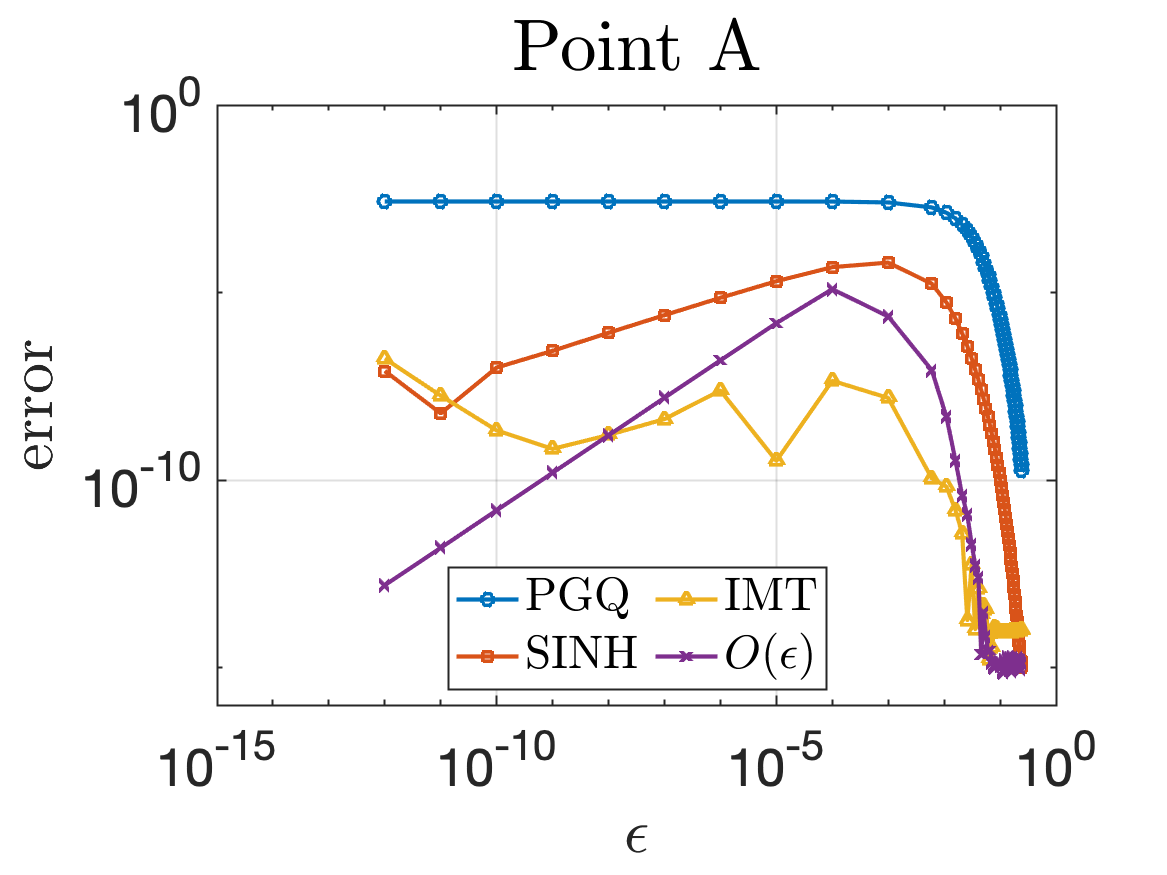}
\includegraphics[width=0.32\linewidth]{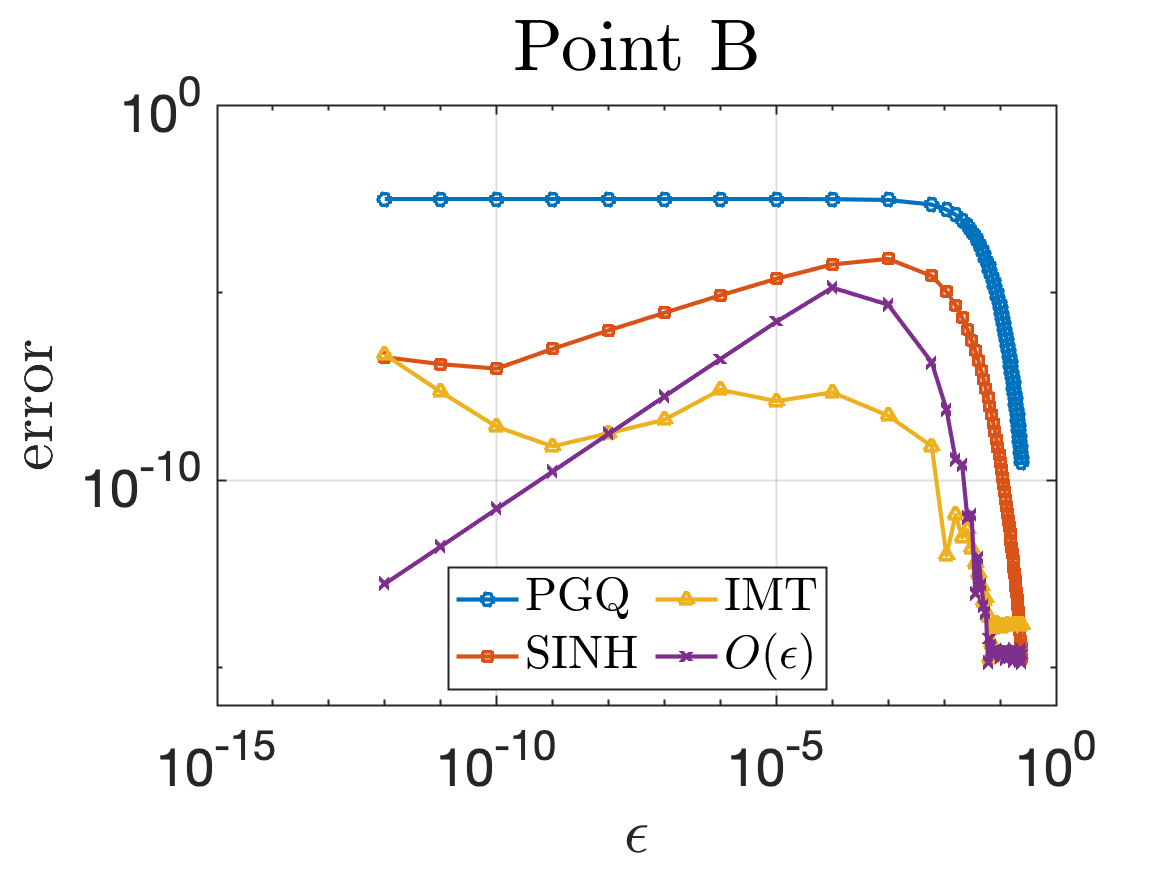}
\includegraphics[width=0.32\linewidth]{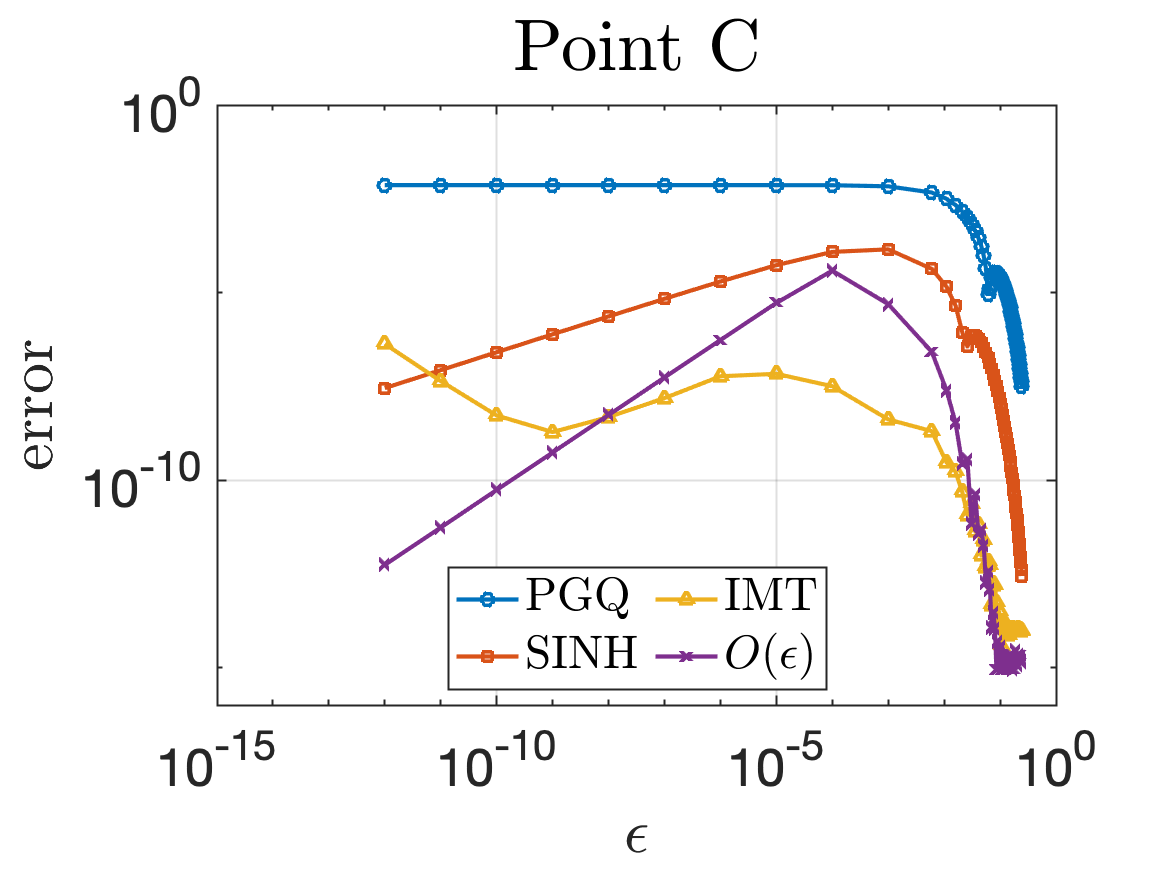}
  \caption{The logarithmic error for {the new numerical method applied to \eqref{eq:repres-int-sub} ($O(\epsilon)$), PGQ method, SINH method, and IMT method with $N = 128$}  at the three zoomed in points, A, B, and C of Fig. \ref{img:mushroom2} when solving the close evaluation problem interior of the mushroom cap domain as a function of $\epsilon$. } 
\label{img:mushroom4}
\end{figure}

 We now show the same analysis within the mushroom cap domain. In Fig. \ref{img:mushroom2}, we present the error interior of the mushroom cap domain when using { the new numerical method with $N = 128$ to approximate \eqref{eq:repres-int-sub}} 
 In Fig. \ref{img:mushroom4} we present the error as a function of $\epsilon$ for all {four} methods {starting} at the three points (A, B, and C) labeled in Fig. \ref{img:mushroom2}. Once again, we observe, as expected, {that our new method has an $O(\epsilon)$ error as  $\epsilon \rightarrow 0^+$.} 
 

{
\subsection{Extension to $O(\epsilon^2)$}
\label{sec:ext} 
}

{
In this section, we present results of when the new numerical method, detailed in Section \ref{sec:numerics}{,} is used to solve both \eqref{eq:repres-int-sub}, the $O(\epsilon)$ method shown above, and \eqref{eq:repres-int-sub-asympt}. In the latter case, we expect the error to be $O(\epsilon^2)$ as $\epsilon \rightarrow 0^+$. 
}

\begin{figure}[h!]
\centering
\includegraphics[width=0.32\linewidth]{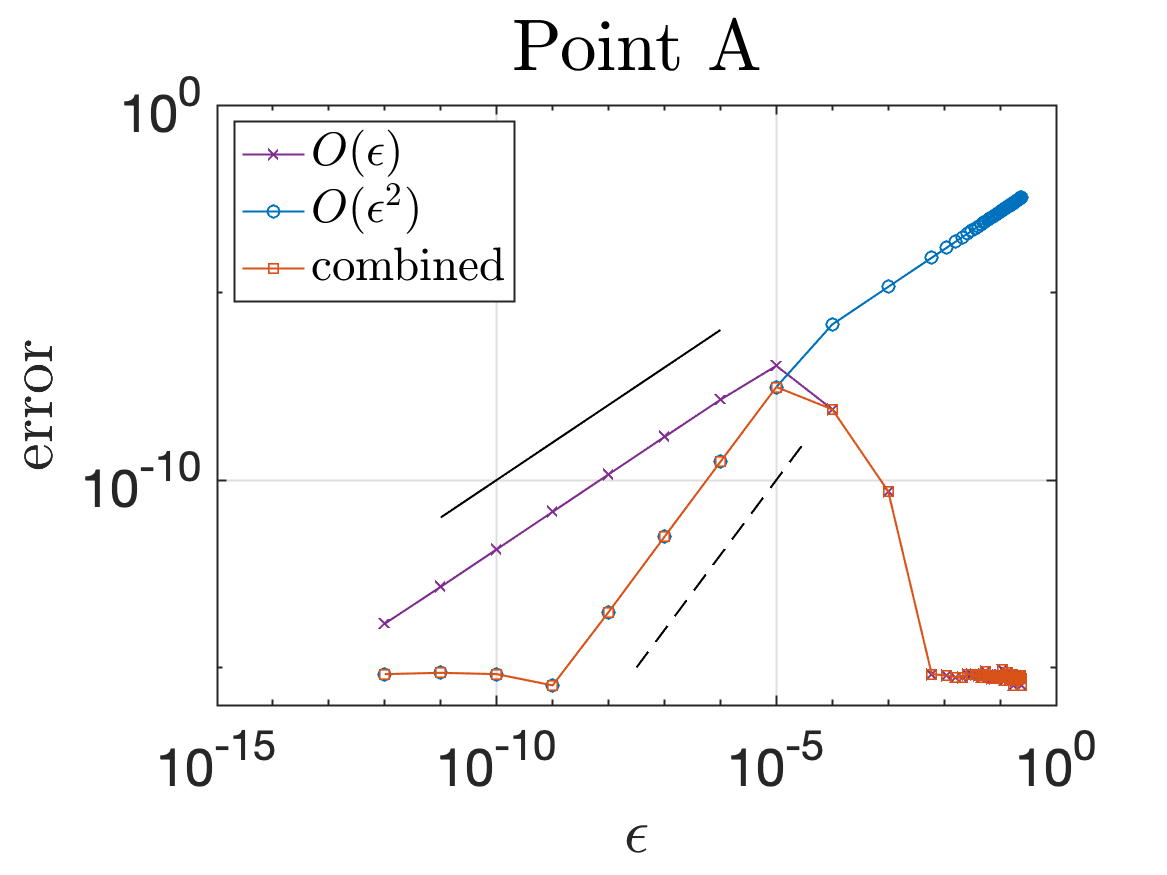}
\includegraphics[width=0.32\linewidth]{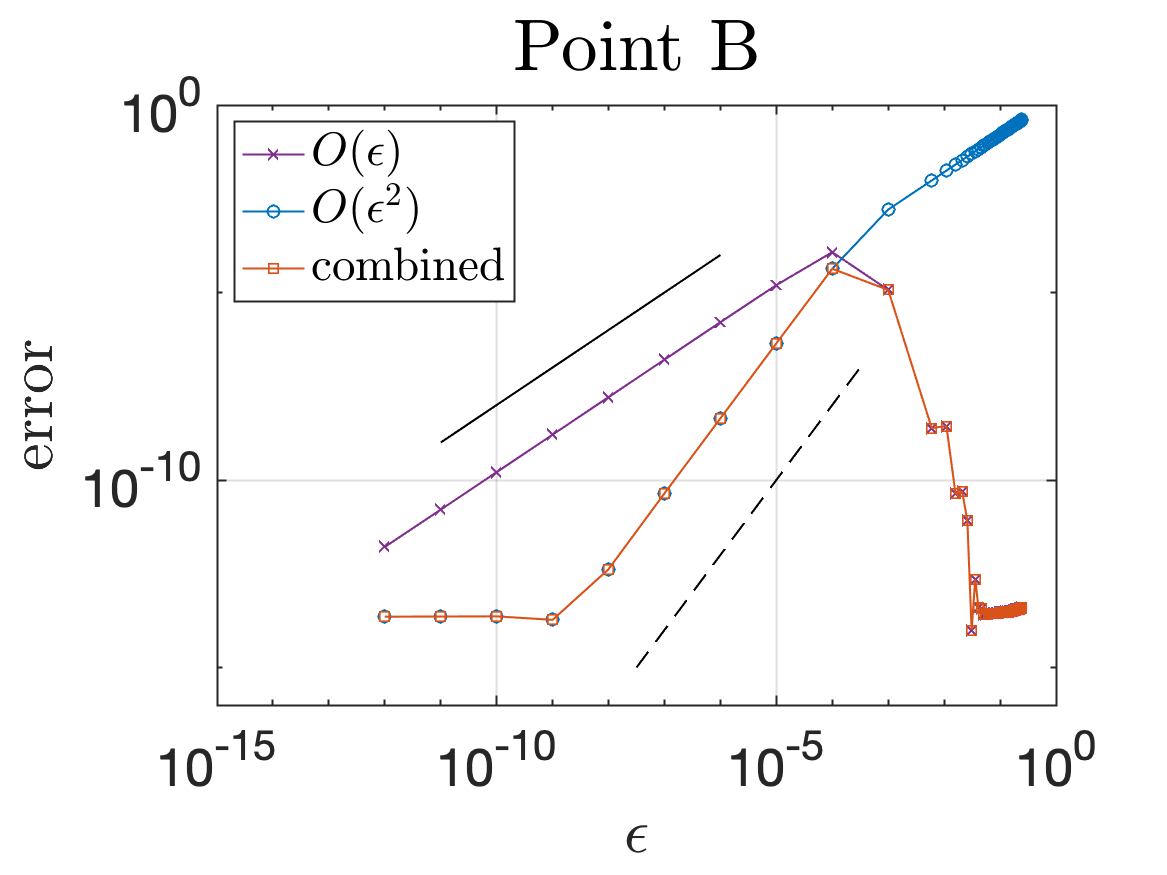}
\includegraphics[width=0.32\linewidth]{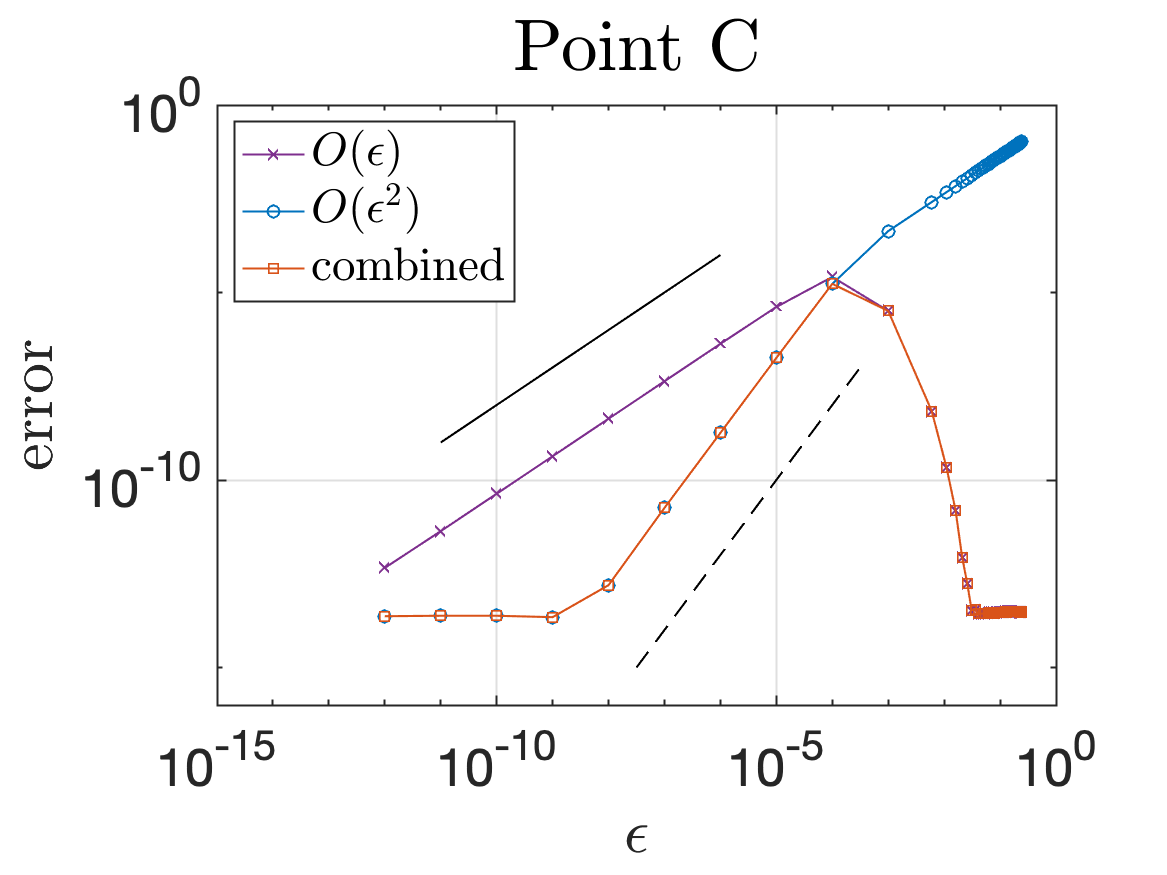}
\caption{{The logarithmic error for the new numerical method applied to \eqref{eq:repres-int-sub} ($O(\epsilon)$) and to \eqref{eq:repres-int-sub-asympt} ($O(\epsilon^2)$) and the combined method with $N = 128$ at the three zoomed in points, A, B, and C of Fig. \ref{img:peanut2} when solving the close evaluation problem interior of the peanut-shaped domain as a function of $\epsilon$.  The solid black line demonstrates $O(\epsilon)$ convergence and the dashed black line demonstrates $O(\epsilon^2)$ convergence.} } 
\label{img:peanut6}
\end{figure}

\begin{figure}[h!]
\centering
\includegraphics[width=0.32\linewidth]{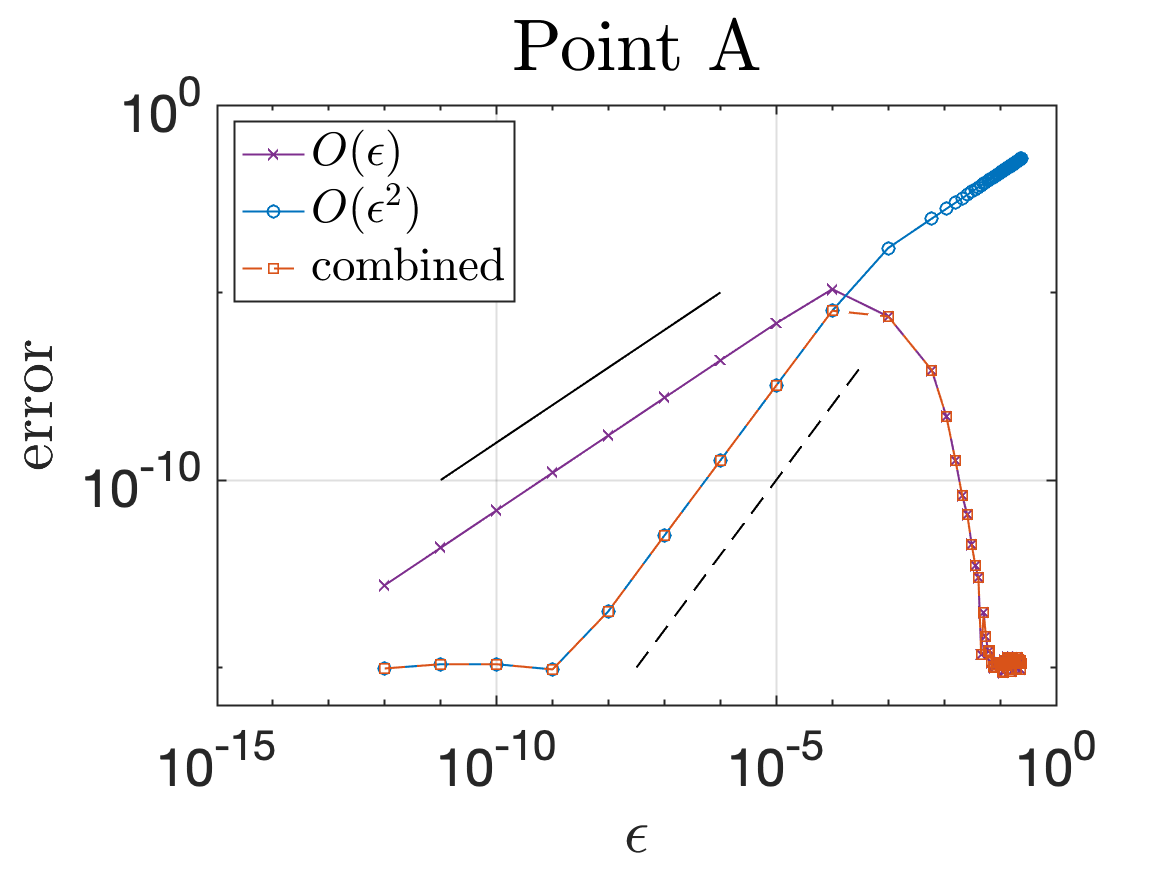}
\includegraphics[width=0.32\linewidth]{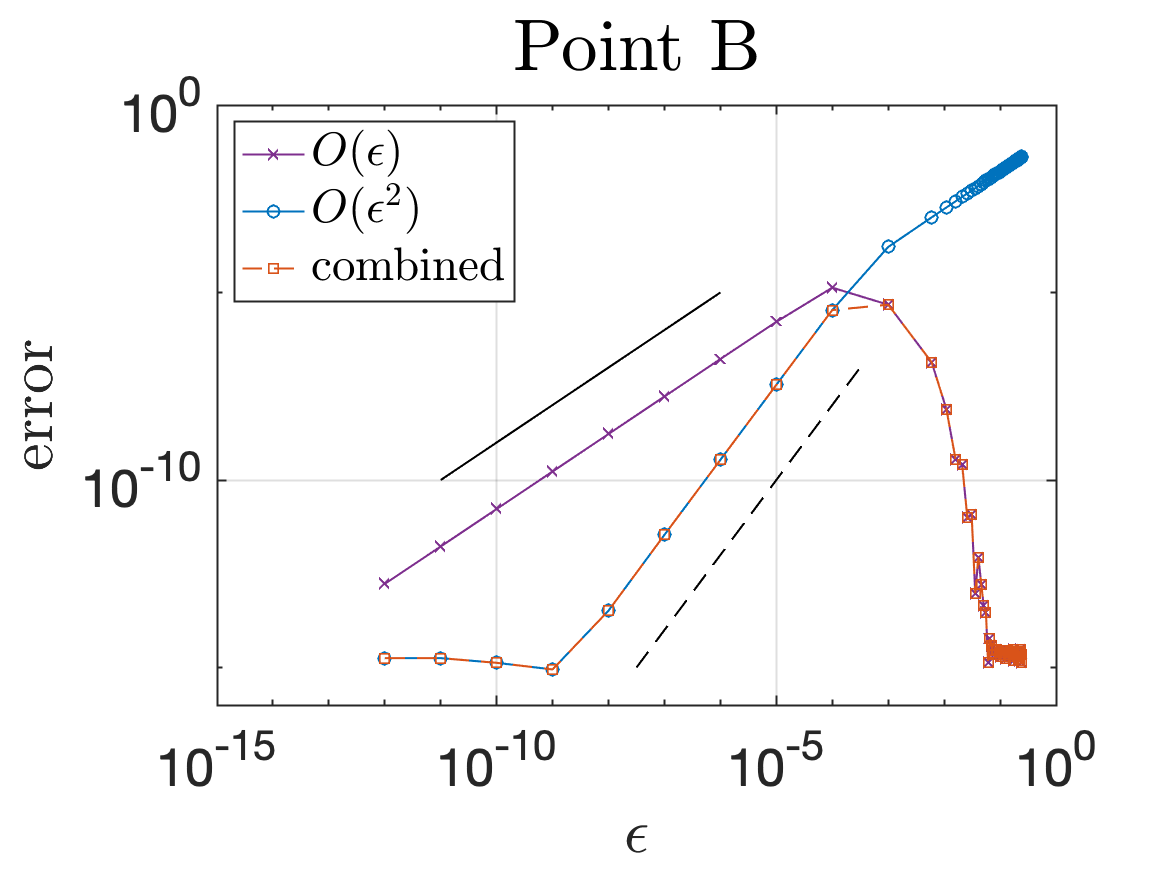}
\includegraphics[width=0.32\linewidth]{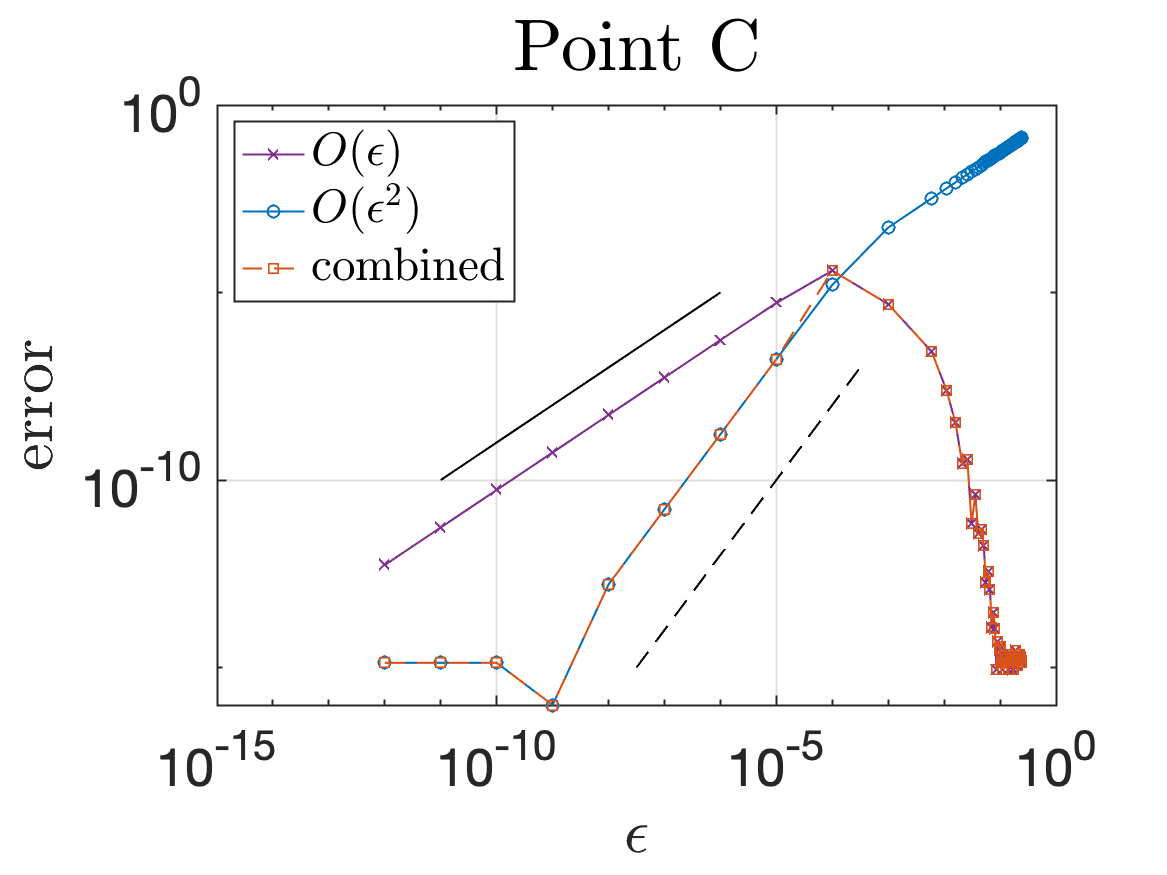}
\caption{{The logarithmic error for the new numerical method applied to \eqref{eq:repres-int-sub} ($O(\epsilon)$) and to \eqref{eq:repres-int-sub-asympt} ($O(\epsilon^2)$) and the combined method with $N = 128$ at the three zoomed in points, A, B, and C of Fig. \ref{img:mushroom2} when solving the close evaluation problem interior of the mushroom cap domain as a function of $\epsilon$. The solid black line demonstrates $O(\epsilon)$ convergence and the dashed black line demonstrates $O(\epsilon^2)$ convergence.}} 
\label{img:mushroom6}
\end{figure}

{
In Fig. \ref{img:peanut6}, we show the logarithmic error for both the $O(\epsilon)$ method and the $O(\epsilon^2)$ method {starting at the} three $y^{\star}$ points (A, B, and C), labeled in Fig. \ref{img:peanut2}, interior to the peanut-shaped domain with $N =128$. Similarly, in Fig. \ref{img:mushroom6}, we show the error of the two methods for the three points labeled in Fig. \ref{img:mushroom2}{,} interior to the mushroom cap domain. Observe that the numerical method applied to \eqref{eq:repres-int-sub-asympt} has an $O(\epsilon^2)$ error as $\epsilon \rightarrow 0^+$ and reaches machine precision for small $\epsilon$. Note that this method is different from the $O(\epsilon)$ method because it depends on an asymptotic expansion of the single-layer potential as $\epsilon \rightarrow 0^+$ and therefore exhibits larger errors as $\epsilon$ increases. 
}

\begin{figure}[h!]
  \centering
   \includegraphics[width=0.32\linewidth]{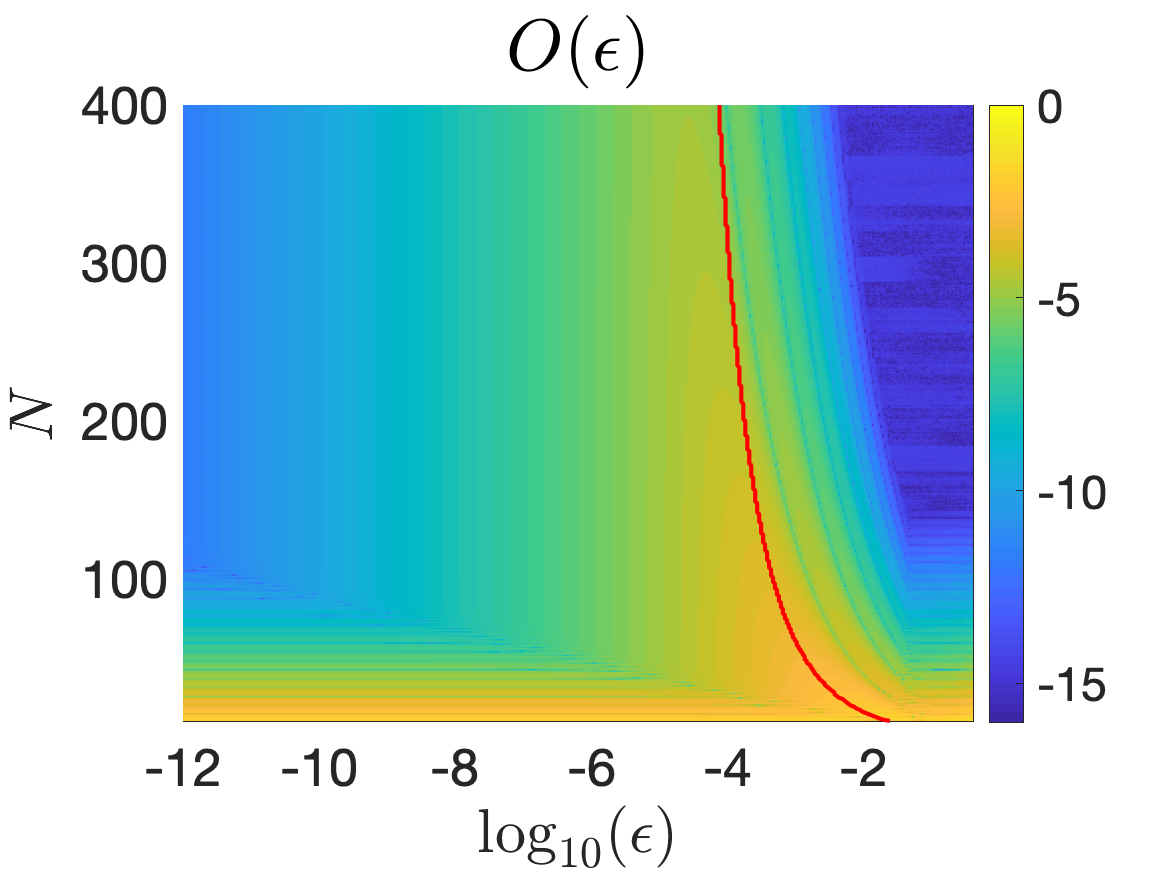}
 \includegraphics[width=0.32\linewidth]{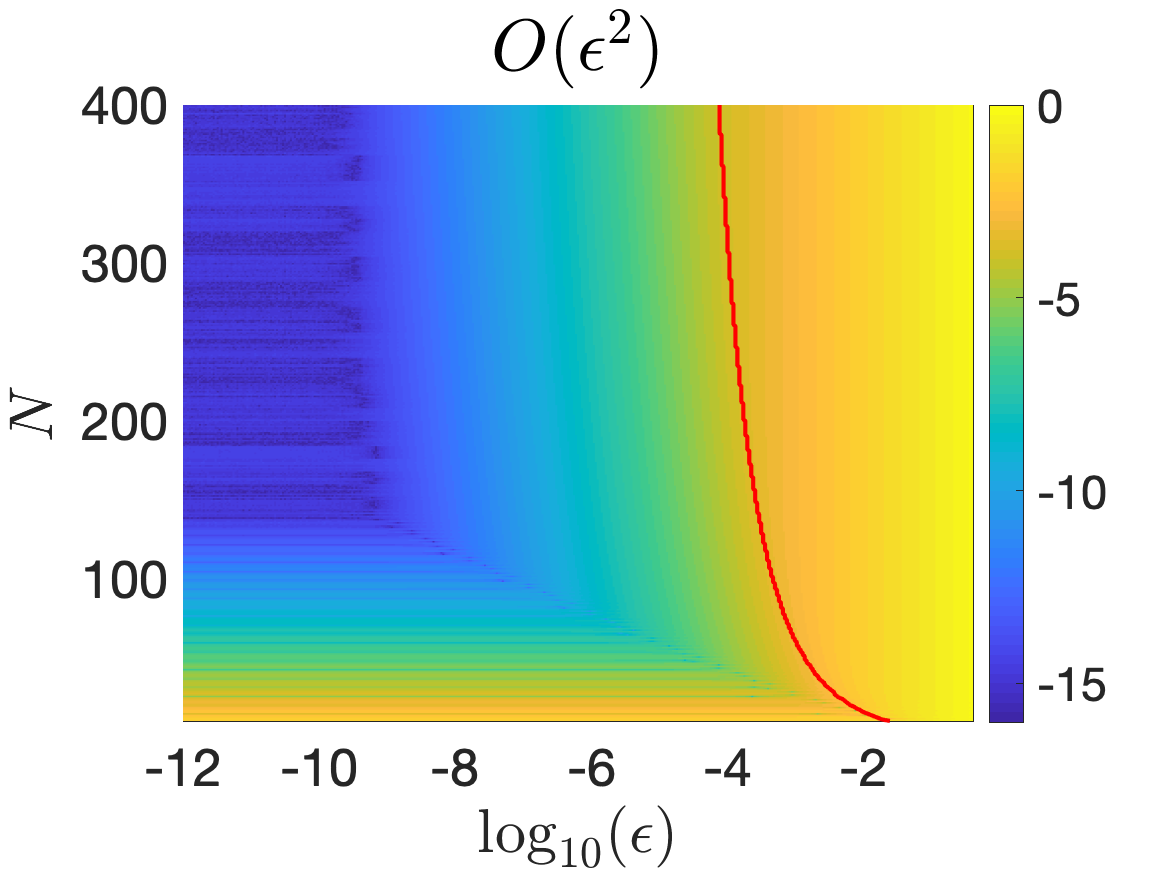}
 \includegraphics[width=0.32\linewidth]{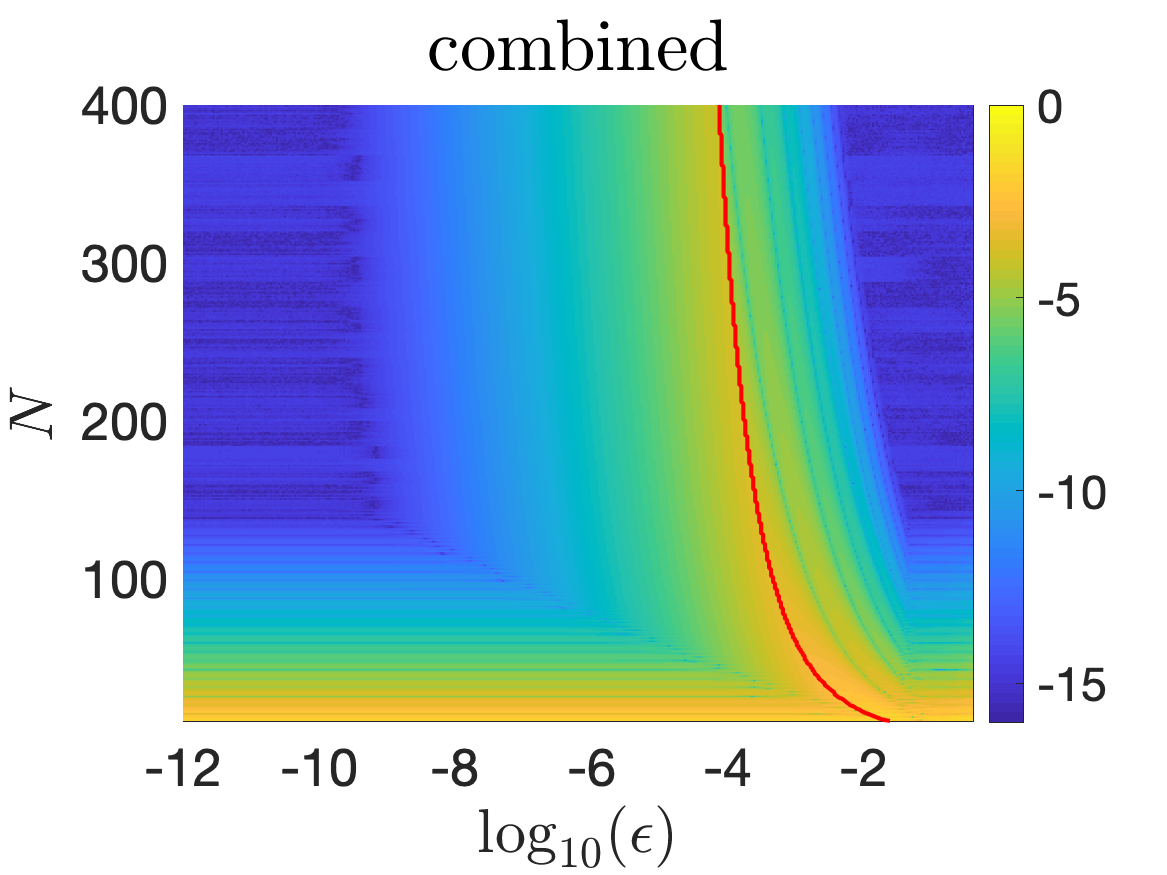}
  \caption{The logarithmic error for the $O(\epsilon)$ method, the $O(\epsilon^2)$ method, and the combined method at the zoomed in point {B} of Fig. \ref{img:peanut2} when solving the close evaluation problem interior of the peanut-shaped domain as a function of $N$ and $\epsilon$. {The red curve on each plot shows where we switch from the $O(\epsilon^2)$ method to the $O(\epsilon)$ method in the combined method.} }
\label{img:peanut7}
\end{figure}

{
Even though the focus of this work is to study the behavior of the close evaluation problem when using a fixed $N$-point quadrature and varying $\epsilon$, we present in Fig. \ref{img:peanut7} the logarithmic error as both $\epsilon$ and $N$ vary for  the $O(\epsilon)$ method and the $O(\epsilon^2)$ method at point B labeled in Fig. \ref{img:peanut2} for the peanut-shaped domain. Similarly, we present  the logarithmic error as both $\epsilon$ and $N$ vary for point B labeled in Fig. \ref{img:mushroom2} for the mushroom cap domain. For both methods, $N$ must be large enough to resolve the domain to not observe $O(1)$ errors. As stated above, for the $O(\epsilon^2)$ method, since we have used an expansion for the single-layer potential, the method exhibits larger errors far from the boundary. For the newly developed methods here, once past a minimum $N$ value ($N \approx 100$ here), the error is constant in $N$ for small values of $\epsilon$. This is expected for the numerical methods developed here based on our asymptotic analysis in $\epsilon$, as presented in Section \ref{sec:local-analysis}. Once again, we observe that the methods are $O(\epsilon)$ and $O(\epsilon^2)$, as expected, and the $O(\epsilon^2)$ method reaches machine precision for small $\epsilon$. 
}

  \begin{figure}[h!]
  \centering
   \includegraphics[width=0.32\linewidth]{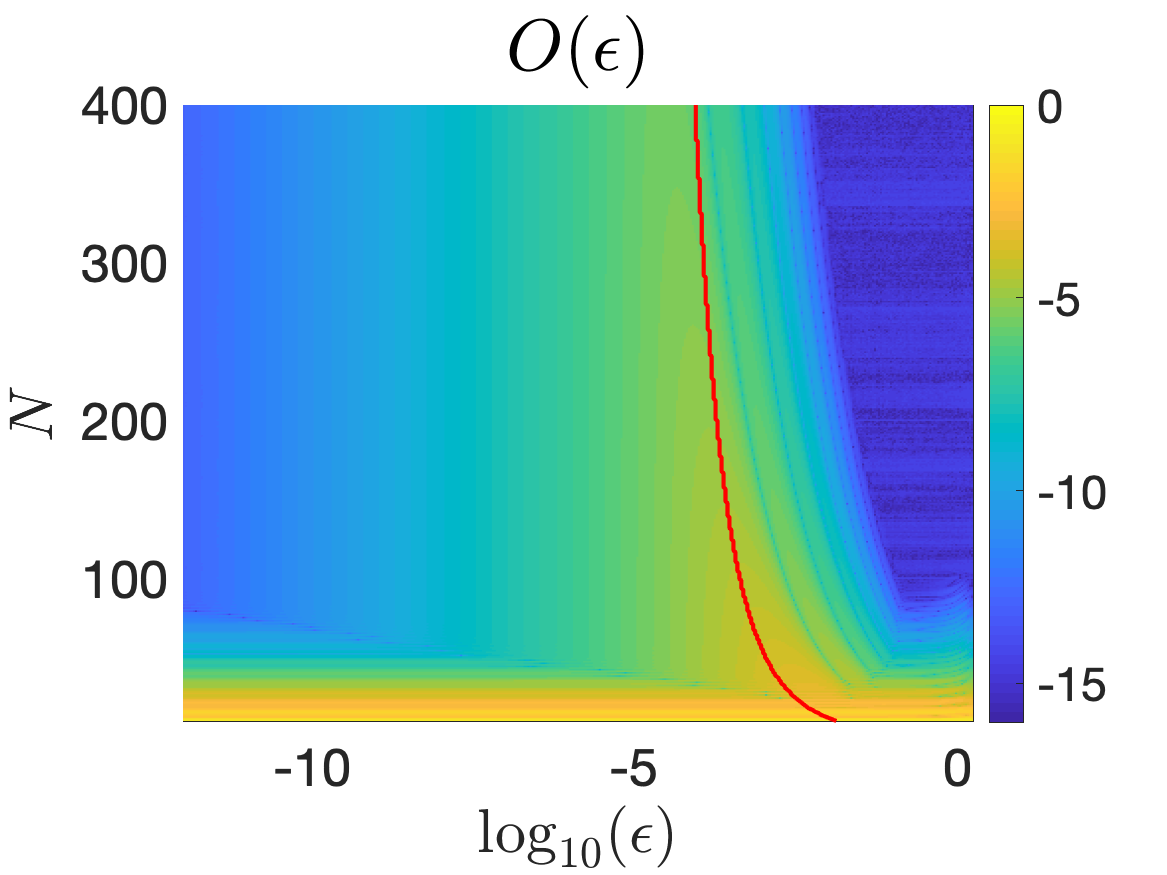}
 \includegraphics[width=0.32\linewidth]{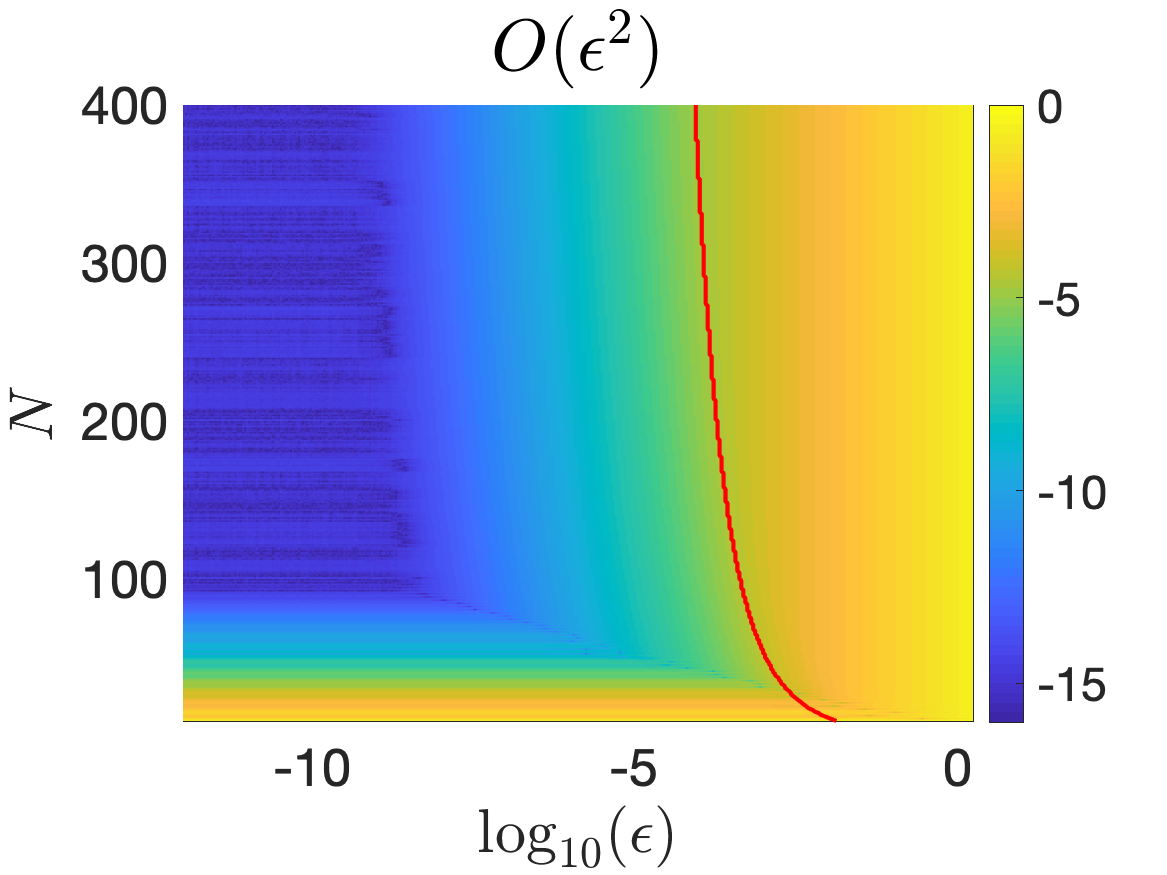}
 \includegraphics[width=0.32\linewidth]{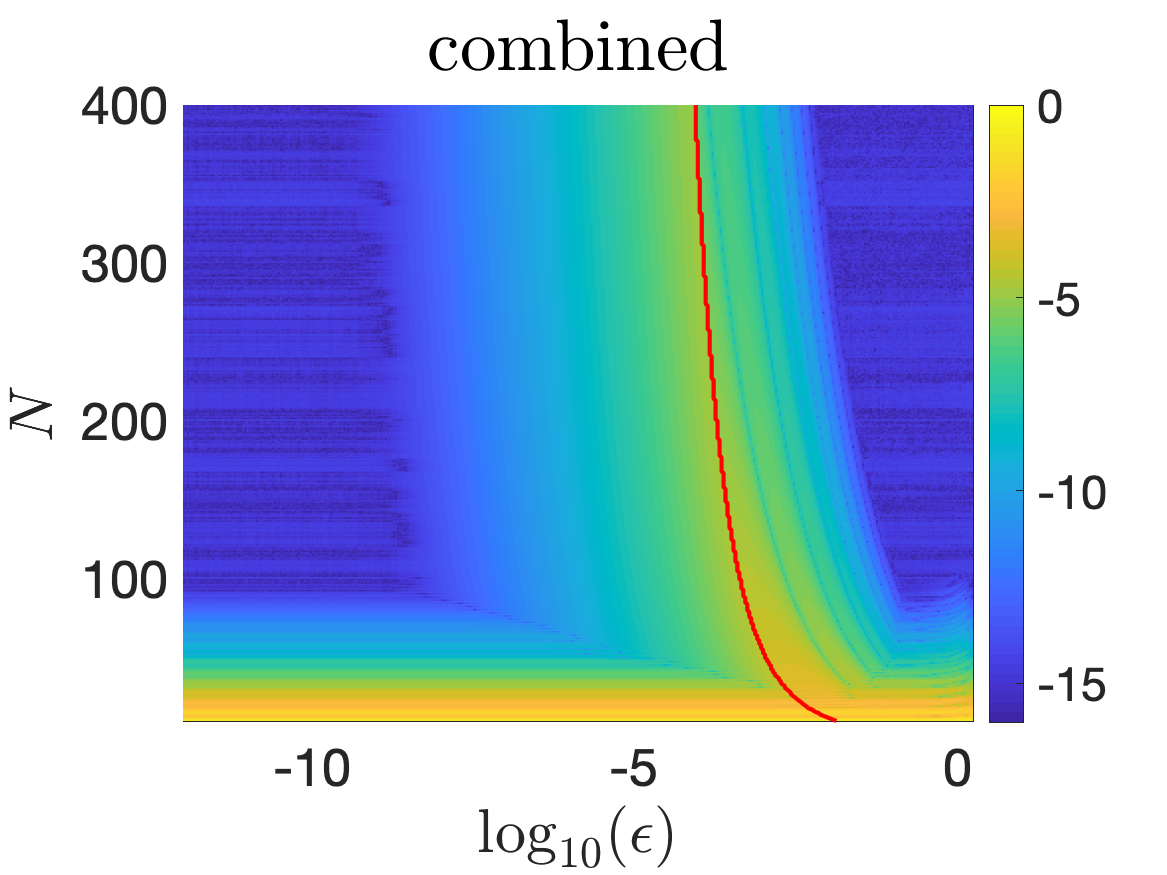}
 
  \caption{The logarithmic error for the $O(\epsilon)$ method, the $O(\epsilon^2)$ method, and the combined method at the zoomed in point {B} of Fig. \ref{img:mushroom2} when solving the close evaluation problem interior of the mushroom cap domain as a function of $N$ and $\epsilon$. {The red curve on each plot shows where we switch from the $O(\epsilon^2)$ method to the $O(\epsilon)$ method in the combined method.} }\label{img:mushroom7}
\end{figure}

{
The optimized numerical method would combine the $O(\epsilon^2)$ method for the smallest values of $\epsilon$ with the $O(\epsilon)$ method for larger values. This combined method is also presented in Figs. \ref{img:peanut6} and \ref{img:mushroom6} and in  Figs.  \ref{img:peanut7} and \ref{img:mushroom7}.  To determine at which value of $\epsilon$ to switch between methods, in particular in general cases when there is not a known exact solution, we extend an idea {we} developed  in two dimensions \cite{carvalho2018asymptotic}  to three dimensions. For a fixed $y^{\star}$ and fixed resolution, $N$, we evaluate Gauss' Law, \eqref{eq:gauss}, for varying $\epsilon$ using the PGQ method and compute the error. This method will suffer from the close evaluation problem and the error approaches $O(1)$ as $\epsilon \rightarrow 0^+$. We choose a tolerance for the error to Gauss' Law, close to where it first approaches $O(1)$. In all cases presented in this paper, this tolerance was chosen to be $0.495$. This error tolerance then corresponds to an $\epsilon$ value for the fixed $y^{\star}$ and fixed resolution, $N$, which is then used as the $\epsilon$ value at which to switch between the two methods. In the plots in Figs. \ref{img:peanut7} and \ref{img:mushroom7}, the red curves give the $\epsilon$ values for each $N$ at which we switch between the $O(\epsilon^2)$ method and the $O(\epsilon)$ method in the combined method. We observe in Figs. \ref{img:peanut6} and \ref{img:mushroom6} and in  Figs.  \ref{img:peanut7} and \ref{img:mushroom7} that the $\epsilon$ value at which we switch between the methods is the ideal value, allowing the combined method to benefit from the best of each of the $O(\epsilon)$  and $O(\epsilon^2)$ methods.
}





\subsection{Effect of Curvature}

\begin{figure}[h!]
\centering
\includegraphics[width=\linewidth]{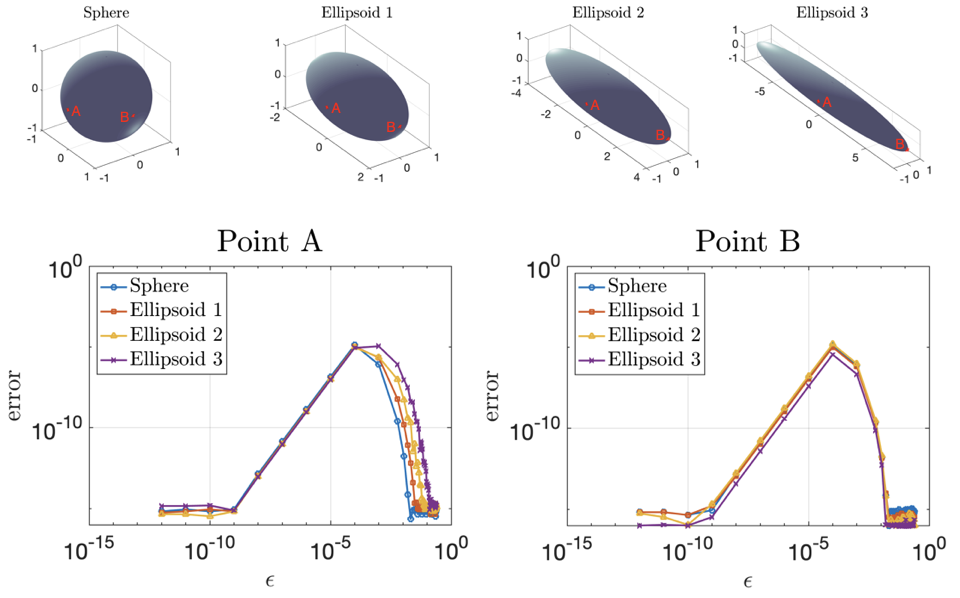}
\caption{{The logarithmic error for the new numerical method with $N = 128$  when solving the close evaluation problem interior to four different domains, a sphere and three ellipsoids (top row) at two points, A: $(-1, 0, 0)$ and B: $(0, b, 0)$ where $b = 1,2,4$ or $8$. }}
\label{fig:ellipsoids}
\end{figure}

{
In the above results presented, we have focused on two domains, the peanut-shaped domain and the mushroom cap domain. It is well known that the curvature plays a role in the error observed in the close evaluation problem \cite{barnett2014evaluation, carvalho2018asymptotic}. Here, we present results to systematically study the robustness of the new method when considering curvature. We test our new method on four different domains, one sphere and three ellipsoids, as shown in Fig. \ref{fig:ellipsoids}, 
\begin{equation}
  y(\theta,\varphi) = ( \sin\theta \cos\varphi, 
 b \sin\theta \sin\varphi, \cos\theta ), \quad \theta \in [0,\pi],
  \quad \varphi \in [-\pi,\pi],
  \label{eq:ellipsoiddomains}
\end{equation}
where $b = 1$ for the sphere and $b = 2,4,8$ for the ellipsoids. In Fig. \ref{fig:ellipsoids}, we present the logarithmic error as a function of $\epsilon$ {starting at} two values of $y^{\star}$ on each boundary, $(-1,0,0)$ (Point A) and  $(0,b,0)$ (Point B) when using the combined method, as discussed above in Section \ref{sec:ext}. 
We observe that our new numerical is robust to variations in curvature and there is no appreciable difference in the errors as the curvature of the boundary increases. 
}

 \subsection{Summary of the results}
 
 For both the peanut-shaped and mushroom cap domains, we have found
 that 
  {the new numerical method we have developed here}
 provides the most accurate evaluation of the
 representation formula for close evaluation points exhibiting an
 error that is $O(\epsilon^2)$ {as  $\epsilon \rightarrow 0^+$}. 
 {The $O(\epsilon)$ method} is the next most
 accurate 
 due to the {increased error from the} 
 single-layer potential. These error results are expected based on the
 analysis provided in Section {\ref{sec:local-analysis}} and the
 extension given in Section \ref{sec:extension}.  These results are
 consistent for two very different three-dimensional domains {and robust when systematically 
 studying the curvature of the boundary}. Thus,
 these results demonstrate that the numerical method developed here
 is accurate and effective for the close evaluation problem.

\section{Conclusions}
\label{sec:conclusions}

We have presented a simple and effective numerical method for the
close evaluation of double- and single-layer potentials in three
dimensions. The close evaluation of these layer potentials are
challenging to numerically compute because they are nearly singular
integrals.  Through a local analysis of the close evaluation of
double- and single-layer potentials about the point at which their
kernels are sharply peaked, we have identified a natural way to
compute these nearly singular integrals. 

Under the simplifying assumption that the boundary {can be parameterized using spherical coordinates}, we work in a rotated coordinate
system in which {the} singular point maps to the north pole of {a} sphere. This rotated coordinate system highlights the axisymmetry of
the kernels for the double- and single-layer potentials at close
evaluation points. We then propose a numerical quadrature rule for the
double- and single-layer potentials written in spherical coordinates. In
this rotated coordinate system, azimuthal integration acts as a
natural averaging operation about the singular point, which enhances its
asymptotic behavior at close evaluation points. The quadrature rule
also integrates over the polar angle directly rather than integrating
over the cosine of the polar angle. By doing so, the numerical
integration includes a factor of sine of the polar angle (the natural
Jacobian for the spherical coordinate system) which is important for
effectively computing these nearly singular integrals. Finally,
because we use an open Gauss-Legendre quadrature rule for the polar
angle, we do not require explicit evaluation of the kernels at the
peaked points.

We have computed results for evaluating the representation
formula interior to a boundary for two different domains and observed the expected behavior.  
{Further, we have shown that our method is robust to variations in curvature through a systematic study.}
 This method can easily be extended to evaluating the representation formula exterior to a boundary. 
 One obtains these results by taking care of the outward normals in the formulas used
throughout this discussion and in Gauss' law.

We have assumed that the boundary is closed, oriented, and analytic in
this paper. However, as long as a portion of the boundary about the
singular point is sufficiently smooth {and can be parameterized using spherical coordinates (in other words integrating over a cap of a sphere)}, these results hold. Thus, this method can be extended
to finite-sized patches covering a surface \cite{bruno2001fast,siegel2018local}.

%

Finally, the analysis shown here can be extended to other problems. In
particular, this method is broadly applicable to weakly singular and
nearly singular integrals over two-dimensional surfaces. These weakly
singular and nearly singular integrals may correspond to solutions of
other elliptic partial differential equations including Helmholtz's
equation and the Stokes' equations. 

\appendix

\section{Rotations on the sphere}
\label{sec:rotation}

We give the explicit rotation formulas over the sphere used throughout
this paper. Consider $v, v^{\star} \in S^{2}$, with $S^{2}$ denoting
the unit sphere. The $(\ihat, \jhat, \khat)$-coordinate system
corresponds to the laboratory reference frame.  We introduce the
parameters $\theta \in [0,\pi]$ and $\varphi \in [-\pi,\pi]$ and write
\begin{equation}
  v = v(\theta,\varphi) = \sin \theta \cos \varphi \, \ihat + \sin
  \theta \sin \varphi \, \jhat + \cos \theta \,
  \khat.
  \label{eq:y-xyz}
\end{equation}
The specific parameters $\theta^{\star}$ and $\varphi^{\star}$ are
defined through the relation
$v^{\star} = v(\theta^{\star},\varphi^{\star})$. We would like to work
in the rotated, $(\hat{e}_{1}, \hat{e}_{2}, \hat{e}_{3})$-coordinate
system in which
\begin{equation}
  \begin{aligned}
    \hat{e}_{1} &= \cos \theta^{\star} \cos \varphi^{\star} \, \ihat +
    \cos \theta^{\star} \sin \varphi^{\star} \, \jhat - \sin
    \theta^{\star} \, \khat,\\
    \hat{e}_{2} &= - \sin \varphi^{\star} \, \ihat + \cos
    \varphi^{\star} \, \jhat,\\
    \hat{e}_{3} &= \sin \theta^{\star} \cos \varphi^{\star} \, \ihat +
    \sin \theta^{\star} \sin \varphi^{\star} \, \jhat + \cos
    \theta^{\star} \, \khat.
  \end{aligned}
  \label{eq:transformations}
\end{equation}
In this rotated system we have $\hat{e}_{3} = v^{\star}$. For this
rotated coordinate system, we introduce the parameters $s \in [0,\pi]$
and $t \in [-\pi,\pi]$ such that
\begin{equation}
  v = v(s,t) = \sin s \cos t \, \hat{e}_{1} + \sin s
  \sin t \, \hat{e}_{2} + \cos s \, \hat{e}_{3}.
  \label{eq:y-uvw}
\end{equation}
It follows that $\hat{e}_{3} = v^{\star} = v(0,\cdot)$. This
corresponds to setting $v^{\star}$ to be the north pole of the rotated
sphere.  By equating \eqref{eq:y-xyz} and \eqref{eq:y-uvw} and
substituting \eqref{eq:transformations} into that result, we obtain
\begin{equation}
  \begin{bmatrix} \sin \theta \cos \varphi\\ \sin \theta \sin \varphi \\
    \cos \theta \end{bmatrix}
  = \begin{bmatrix} \cos \theta^{\star} \cos \varphi^{\star} & - \sin
    \varphi^{\star} & \sin \theta^{\star} \cos \varphi^{\star}\\
    \cos \theta^{\star} \sin \varphi^{\star} & \cos \varphi^{\star} & \sin
    \theta^{\star} \sin \varphi^{\star} \\
    - \sin \theta^{\star} & 0 & \cos \theta^{\star}
  \end{bmatrix}
  \begin{bmatrix} \sin s \cos t\\ \sin s \sin t\\
    \cos s \end{bmatrix}.
  \label{eq:y-identity}
\end{equation}
Let us rewrite \eqref{eq:y-identity} compactly as
$v(\theta,\varphi) = R(\theta^{\star},\varphi^{\star}) v(s,t)$ with
$R(\theta^{\star},\varphi^{\star})$ denoting the $3 \times 3$ matrix
given above. It is a rotation matrix. Hence, it is orthogonal.

We now seek to write $\theta = \theta(s,t)$ and
$\varphi = \varphi(s,t)$. To do so, we introduce
\begin{align}
  \xi(s,t;\theta^{\star},\varphi^{\star}) 
  &= \cos \theta^{\star} \cos \varphi^{\star} \sin s \cos t -
    \sin \varphi^{\star} \sin s \sin t + \sin \theta^{\star} \cos
    \varphi^{\star} \cos s,  \label{eq:angle-rotations1} \\
  \eta(s,t;\theta^{\star},\varphi^{\star})
  &= \cos \theta^{\star} \sin \varphi^{\star} \sin s \cos t + \cos
    \varphi^{\star} \sin s \sin t + \sin \theta^{\star} \sin
    \varphi^{\star} \cos s,\\
  \zeta(s,t;\theta^{\star},\varphi^{\star})
  &= - \sin \theta^{\star} \sin s \cos t + \cos \theta^{\star} \cos s.
    \label{eq:angle-rotations3}
\end{align}
From \eqref{eq:angle-rotations1} - \eqref{eq:angle-rotations3}, we
find that
\begin{equation}
  \theta = \arctan \left( \frac{\sqrt{\xi^{2} + \eta^{2}}}{\zeta}
  \right),
  \label{eq:theta-st}
\end{equation}
and
\begin{equation}
  \varphi = \arctan\left( \frac{\eta}{\xi} \right).
  \label{eq:varphi-st}
\end{equation}
These results give the formulas needed to evaluate
$\theta = \theta(s_i,t_j)$ and $\varphi = \varphi(s_i,t_j)$ that are used in 
Section{s \ref{sec:prior} and \ref{sec:numerics}}. 
{This $\theta$ and $\phi$ are used to evaluate $y(s_i,t_j) = y(\theta(s_i,t_j),  \varphi(s_i,t_j)) $ using \eqref{eq:generaldomain} or \eqref{eq:ellipsoiddomains} in the quadrature rules  \eqref{eq:genquadrature} and \eqref{eq:quadrature}.}


\section*{Acknowledgments}

This research was supported by National Science Foundation Grant:
DMS-1819052. A. D. Kim is also supported by the Air Force Office of
Scientific Research Grants: FA9550-17-1-0238 and FA9550-18-1-0519.
S. Khatri is also supported by the National Science Foundation Grant: PHY-1505061. R. Cortez is partially supported by National Science Foundation Grant: DMS-1043626. 

\end{document}

%% file: sketch.pdf_tex
\begingroup%
  \makeatletter%
  \providecommand\color[2][]{%
    \errmessage{(Inkscape) Color is used for the text in Inkscape, but the package 'color.sty' is not loaded}%
    \renewcommand\color[2][]{}%
  }%
  \providecommand\transparent[1]{%
    \errmessage{(Inkscape) Transparency is used (non-zero) for the text in Inkscape, but the package 'transparent.sty' is not loaded}%
    \renewcommand\transparent[1]{}%
  }%
  \providecommand\rotatebox[2]{#2}%
  \ifx\svgwidth\undefined%
    \setlength{\unitlength}{228.53649414bp}%
    \ifx\svgscale\undefined%
      \relax%
    \else%
      \setlength{\unitlength}{\unitlength * \real{\svgscale}}%
    \fi%
  \else%
    \setlength{\unitlength}{\svgwidth}%
  \fi%
  \global\let\svgwidth\undefined%
  \global\let\svgscale\undefined%
  \makeatother%
  \begin{picture}(1,0.77127729)%
    \put(0,0){\includegraphics[width=\unitlength,page=1]{sketch.pdf}}%
    \put(0.39869752,0.52703723){\color[rgb]{0,0,0}\makebox(0,0)[lt]{\begin{minipage}{0.32880026\unitlength}\raggedright $y^\star$\end{minipage}}}%
    \put(0.49859344,0.04356715){\color[rgb]{0,0,0}\makebox(0,0)[lt]{\begin{minipage}{0.19377958\unitlength}\raggedright $x$\end{minipage}}}%
    \put(0,0){\includegraphics[width=\unitlength,page=2]{sketch.pdf}}%
    \put(0.28743047,0.77957576){\color[rgb]{0,0,0}\makebox(0,0)[lt]{\begin{minipage}{0.34380252\unitlength}\raggedright $n^\star$\end{minipage}}}%
    \put(0.80469103,0.42924306){\color[rgb]{0,0,0}\makebox(0,0)[lt]{\begin{minipage}{0.23878651\unitlength}\raggedright $B$\end{minipage}}}%
    \put(0.19798663,0.27787222){\color[rgb]{0,0,0}\makebox(0,0)[lt]{\begin{minipage}{0.19377961\unitlength}\raggedright $D$\end{minipage}}}%
    \put(0.35497208,0.18604407){\color[rgb]{0,0,0}\makebox(0,0)[lt]{\begin{minipage}{0.60820479\unitlength}\raggedright $\epsilon \ell$\end{minipage}}}%
    \put(0,0){\includegraphics[width=\unitlength,page=3]{sketch.pdf}}%
    \put(-0.59407807,0.4779925){\color[rgb]{0,0,0}\makebox(0,0)[lb]{\smash{}}}%
  \end{picture}%
\endgroup%

%% file: KKCC3Drevision.bbl
\begin{thebibliography}{10}
\expandafter\ifx\csname url\endcsname\relax
  \def\url#1{\texttt{#1}}\fi
\expandafter\ifx\csname urlprefix\endcsname\relax\def\urlprefix{URL }\fi

\bibitem{af2016fast}
L.~af~Klinteberg, A.-K. Tornberg, A fast integral equation method for solid
  particles in viscous flow using quadrature by expansion, J. Comput. Phys. 326
  (2016) 420--445.

\bibitem{af2017error}
L.~af~Klinteberg, A.-K. Tornberg, Error estimation for quadrature by expansion
  in layer potential evaluation, Adv. Comput. Math. 43~(1) (2017) 195--234.

\bibitem{atkinson1982numerical}
K.~E. Atkinson, Numerical integration on the sphere, ANZIAM J. 23~(3) (1982)
  332--347.

\bibitem{atkinson1982laplace}
K.~E. Atkinson, The numerical solution {L}aplace's equation in three
  dimensions, SIAM J. Numer. Anal. 19~(2) (1982) 263--274.

\bibitem{atkinson1985algorithm}
K.~E. Atkinson, Algorithm 629: An integral equation program for {L}aplace's
  equation in three dimensions, ACM Trans. Math. Softw. 11~(2) (1985) 85--96.

\bibitem{atkinson1997numerical}
K.~E. Atkinson, The Numerical Solution of Integral Equations of the Second
  Kind, Cambridge University Press, 1997.

\bibitem{barnett2014evaluation}
A.~H. Barnett, Evaluation of layer potentials close to the boundary for
  {L}aplace and {H}elmholtz problems on analytic planar domains, SIAM J. Sci.
  Comput. 36~(2) (2014) A427--A451.

\bibitem{beale2001method}
J.~T. Beale, M.-C. Lai, A method for computing nearly singular integrals, SIAM
  J. Numer. Anal. 38~(6) (2001) 1902--1925.

\bibitem{beale2016simple}
J.~T. Beale, W.~Ying, J.~R. Wilson, A simple method for computing singular or
  nearly singular integrals on closed surfaces, Commun. Comput. Phys. 20~(3)
  (2016) 733--753.

\bibitem{bremer2010nonlinear}
J.~Bremer, Z.~Gimbutas, V.~Rokhlin, A nonlinear optimization procedure for
  generalized gaussian quadratures, SIAM J. Sci. Comput. 32~(4) (2010)
  1761--1788.

\bibitem{bruno2001fast}
O.~P. Bruno, L.~A. Kunyansky, A fast, high-order algorithm for the solution of
  surface scattering problems: basic implementation, tests, and applications,
  J. Comput. Phys. 169~(1) (2001) 80--110.

\bibitem{carvalho2018asymptotic}
C.~Carvalho, S.~Khatri, A.~D. Kim, Asymptotic analysis for close evaluation of
  layer potentials, J. Comput. Phys. 355 (2018) 327--341.

\bibitem{ckk2018asymp}
C.~Carvalho, S.~Khatri, A.~D. Kim, Asymptotic approximations for the close
  evaluation of double-layer potentials, SIAM J. Sci. Comput. 42~(1) (2020)
  A504--A533.

\bibitem{delves1988computational}
L.~M. Delves, J.~L. Mohamed, Computational Methods for Integral Equations,
  Cambridge University Press, 1988.

\bibitem{epstein2013convergence}
C.~L. Epstein, L.~Greengard, A.~Kl{\"o}ckner, On the convergence of local
  expansions of layer potentials, SIAM J. Numer. Anal. 51~(5) (2013)
  2660--2679.

\bibitem{folland1995introduction}
G.~B. Folland, \textcolor{black}{Introduction to partial differential
  equations}, vol. \textcolor{black}{102}, \textcolor{black}{Princeton
  university press}, \textcolor{black}{1995}.

\bibitem{ganesh2004high}
M.~Ganesh, I.~Graham, A high-order algorithm for obstacle scattering in three
  dimensions, J. Comput. Phys. 198~(1) (2004) 211--242.

\bibitem{gimbutas2013fast}
Z.~Gimbutas, S.~Veerapaneni, A fast algorithm for spherical grid rotations and
  its application to singular quadrature, SIAM J. Sci. Comput. 35~(6) (2013)
  A2738--A2751.

\bibitem{graham2002fully}
I.~G. Graham, I.~H. Sloan, \textcolor{black}{Fully discrete spectral boundary
  integral methods for Helmholtz problems on smooth closed surfaces in
  $\mathbb{R}^3$}, \textcolor{black}{Numerische Mathematik}
  \textcolor{black}{92}~(\textcolor{black}{2}) (\textcolor{black}{2002})
  \textcolor{black}{289--323}.

\bibitem{FMM}
L.~Greengard, V.~Rokhlin, A new version of the {F}ast {M}ultipole {M}ethod for
  the {L}aplace equation in three dimensions, Acta Numerica 6 (1997) 229--269.

\bibitem{guenther1996partial}
R.~B. Guenther, J.~W. Lee, Partial Differential Equations of Mathematical
  Physics and Integral Equations, Dover Publications, 1996.

\bibitem{helsing2008evaluation}
J.~Helsing, R.~Ojala, On the evaluation of layer potentials close to their
  sources, J. Comput. Phys. 227~(5) (2008) 2899--2921.

\bibitem{iri1987certain}
M.~Iri, S.~Moriguti, Y.~Takasawa, \textcolor{black}{On a certain quadrature
  formula}, \textcolor{black}{J. Comput. Appl. Math}
  \textcolor{black}{17}~(\textcolor{black}{1-2}) (\textcolor{black}{1987})
  \textcolor{black}{3--20}.

\bibitem{johnston2005sinh}
P.~R. Johnston, D.~Elliott, \textcolor{black}{A sinh transformation for
  evaluating nearly singular boundary element integrals},
  \textcolor{black}{Int. J. Numer. Meth. Eng.}
  \textcolor{black}{62}~(\textcolor{black}{4}) (\textcolor{black}{2005})
  \textcolor{black}{564--578}.

\bibitem{klockner2013quadrature}
A.~Kl{\"o}ckner, A.~Barnett, L.~Greengard, M.~O'Neil, Quadrature by expansion:
  A new method for the evaluation of layer potentials, J. Comput. Phys. 252
  (2013) 332--349.

\bibitem{rachh2017fast}
M.~Rachh, A.~Kl{\"o}ckner, M.~O'Neil, Fast algorithms for quadrature by
  expansion i: Globally valid expansions, J. Comput. Phys. 345 (2017) 706--731.

\bibitem{robinson1981algorithm}
I.~Robinson, E.~De~Doncker, Algorithm 45. {A}utomatic computation of improper
  integrals over a bounded or unbounded planar region, Computing 27~(3) (1981)
  253--284.

\bibitem{schwab1999extraction}
C.~Schwab, W.~Wendland, On the extraction technique in boundary integral
  equations, Math. Comput. 68~(225) (1999) 91--122.

\bibitem{siegel2018local}
M.~Siegel, A.-K. Tornberg, A local target specific quadrature by expansion
  method for evaluation of layer potentials in 3d, J. Comput. Phys. 364 (2018)
  365--392.

\bibitem{wala20183dqbx}
M.~Wala, A.~Kl{\"o}ckner, A fast algorithm for quadrature by expansion in three
  dimensions, \textcolor{black}{J. Comput. Phys.} \textcolor{black}{388}
  (\textcolor{black}{2019}) \textcolor{black}{655--689}.

\end{thebibliography}
